\documentclass[final]{siamart1116}

\usepackage{lipsum}
\usepackage{amsfonts}
\usepackage{graphicx}
\usepackage{epstopdf}
\usepackage{algorithmic}
\usepackage{enumitem}
\ifpdf
  \DeclareGraphicsExtensions{.eps,.pdf,.png,.jpg}
\else
  \DeclareGraphicsExtensions{.eps}
\fi
\usepackage{graphicx} 
\usepackage{float}

\usepackage{algorithm}
\usepackage{algorithmic}

\usepackage{hyperref}

\usepackage{amsmath,amssymb}
\usepackage{pgfplots}
\usepackage{color}
\usepackage{needspace}

\usepackage{flushend}
\usepackage{multirow}
\usepackage{rotating}
\usepackage{appendix}
\usepackage{tikz,pgfplots}
\usepackage{subcaption}
\usepackage{bbm}

\usepackage{standalone} 

\newif\ifcompiletikz
\compiletikztrue

\input{mysymbol.sty}

\newtheorem{assumption}{\hspace{0pt}\bf Assumption}

\numberwithin{theorem}{section}

\newcommand{\TheTitle}{A Newton-Based Method for Nonconvex Optimization with Fast Evasion of Saddle Points}
\newcommand{\shortTitle}{A Newton-Based Method for Nonconvex Optimization} 
\newcommand{\TheAuthors}{Santiago Paternain, Aryan Mokhtari, and Alejandro Ribeiro}

\headers{\shortTitle}{\TheAuthors}

\title{{\TheTitle}\thanks{Submitted to the editors on September 30, 2017. Revised on July 20, 2018.
\funding{Work supported by the ARL DCIST CRA W911NF-17-2-0181}}}

\author{
  Santiago Paternain\thanks{Department of Electrical and Systems Engineering, University of Pennsylvania, Philadelphia, PA 19104
    (\email{spater@seas.upenn.edu}, \email{aribeiro@seas.upenn.edu}).}
  \and
  Aryan Mokhtari\thanks{Laboratory for Information and Decision Systems, Massachusetts Institute of Technology, Cambridge, MA 02139 (\email{aryanm@mit.edu}) }
  \and
  Alejandro Ribeiro\footnotemark[2]
  }

\usepackage{amsopn}

\ifpdf
\hypersetup{
  pdftitle={\TheTitle},
  pdfauthor={\TheAuthors}
}
\fi

\begin{document}
\maketitle

%
\begin{abstract}	 
Machine learning problems such as neural network training, tensor decomposition, and matrix factorization, require local minimization of a nonconvex function. This local minimization is challenged by the presence of saddle points, of which there can be many and from which descent methods may take an inordinately large number of iterations to escape. This paper presents a second-order method that modifies the update of Newton's method by replacing the negative eigenvalues of the Hessian by their absolute values and uses a truncated version of the resulting matrix to account for the objective's curvature. The method is shown to escape saddles in at most $1 + \log_{3/2} (\delta/2\varepsilon)$ iterations where $\varepsilon$ is the target optimality and $\delta$ characterizes a point sufficiently far away from the saddle. This base of this exponential escape is $3/2$ independently of problem constants. Adding classical properties of Newton's method, the paper proves convergence to a local minimum with probability $1-p$ in $O\left(\log(1/p)) + O(\log(1/\varepsilon)\right)$ iterations.
\end{abstract}




%
\section{Introduction}
Although it is generally accepted that the distinction between functions that are easy and difficult to minimize is their convexity, a more accurate statement is that the distinction lies on the ability to use local descent methods. A convex function is easy to minimize because a minimum can be found by following local descent directions, but this is not possible for nonconvex functions. This is unfortunate because many interesting problems in machine learning can be reduced to the minimization of nonconvex functions \cite{pardalos1994nonconvex}. Despite this general complexity, some recent results have shown that for a large class of nonconvex problems such as dictionary learning \cite{sun2017complete}, tensor decomposition \cite{ge2015escaping}, matrix completion \cite{ge2016matrix}, and training of some specific forms of neural networks \cite{kawaguchi2016deep}, all local minimizers are global minima. This reduces the problem of finding {\it the} global optimum to the problem of finding {\it a} local minimum which can be accomplished with local descent methods. 

Conceptually, finding a local minimum of a nonconvex function is not more difficult than finding the minimum of a convex function. It is true that the former can have saddle points that are attractors of gradient fields for some initial conditions \cite[Section 1.2.3]{nesterov2013introductory}. However, since these initial conditions lie in a low dimensional manifold, gradient descent can be shown to converge almost surely to a local minimum if the initial condition is assumed randomly chosen \cite{lee2016gradient, panageas2016gradient}, or if noise is added to gradient descent steps \cite{pemantle1990nonconvergence}. These fundamental facts notwithstanding, practical implementations show that finding a local minimum of a nonconvex function is much more challenging than finding the minimum of a convex function. This happens because the performance of first order methods is degraded by ill conditioning which in the case of nonconvex functions implies that it may take a very large number of iterations to escape from a saddle point \cite{dauphin2014identifying, choromanska2015loss}. Indeed, it can be argued that it is saddle-points and not local minima that provide a fundamental impediment to rapid high dimensional non-convex optimization \cite{dauphin2014identifying, baldi1989neural, saxe2013learning, saxe2013exact}.

In this paper we propose the nonconvex Newton (NCN) method to accelerate the speed of escaping saddles. NCN uses a descent direction analogous to the Newton step except that we use the {\it Positive definite Truncated (PT)-inverse} of the Hessian in lieu of the regular inverse of the Hessian (Definition \ref{def_trucho_inverse}). The PT-inverse has the same eigenvector basis of the regular inverse but its eigenvalues differ in that: (i) All negative eigenvalues are replaced by their absolute values. (ii) Small eigenvalues are replaced by a constant. The idea of using the absolute value of the eigenvalues of the Hessian in nonconvex optimization was first proposed in \cite[Chapters 4 and 7]{nocedal2006numerical} and then in \cite{murray2010newton,pascanu2014saddle}. These properties ensure that the value of the function is reduced at each iteration with an appropriate selection of the step size. Our main contribution is to show that NCN can escape any saddle point with eigenvalues bounded away from zero at an exponential rate which can be further shown to have a base of 3/2 independently of the function's properties in a neighborhood of the saddle. Specifically, we show the following result:

\medskip\begin{itemize}
\item[(i)] Consider an arbitrary $\varepsilon>0$ and the region around a saddle at which the objective gradient is smaller than $\varepsilon$. There exists a subset of this region so that NCN iterations result in the norm of the gradient growing from $\varepsilon$ to $\delta$ at an exponential rate with base $3/2$. The number of NCN iterations required for the gradient to progress from $\varepsilon$ to $\delta$ is therefore not larger than $1 + \log_{3/2} (\delta/2\varepsilon)$; see Theorem \ref{thm_steps_escape}.
\end{itemize}

\medskip \noindent We emphasize that the base of escape 3/2 is independent of the function's properties asides from the requirement to have non-degenerate saddles. The constant $\delta$ depends on Lipschitz constants and Hessian eigenvalue bounds.  

As stated in (i) the base 3/2 for exponential escape does not hold for all points close to the saddle but in a specific subset at which the gradient norm is smaller than $\varepsilon$. It is impossible to show that NCN iterates stay within this region as they approach the saddle, but we show that it is possible to add noise to NCN iterates to quickly enter into this subset with overwhelming probability. Specifically, we show that: 

\medskip\begin{itemize}
\item[(ii)] By adding gaussian noise with standard deviation proportional to $\varepsilon$ when the norm of the gradient of the function is smaller than $\varepsilon$,  the region in which the base of the exponential escape of NCN is $3/2$ is visited by the iterates with probability $1-p$ in $O(\log(1/p))$ iterations. Once this region is visited once, result (i) holds and we escape the saddle in not more than $1 + \log_{3/2} (\delta/2\varepsilon)$ iterations; see Proposition \ref{prop_probabilistic}.
\end{itemize}

\medskip \noindent  Combined with other standard properties of classical Newton's method, results (i) and (ii) imply convergence to a local minimum with probability $1-p$ in a number of iterations that is of order $O(\log(1/\varepsilon))$ with respect to the target accuracy and of order $O\left(\log(1/p)\right)$ with respect to the desired probability (Theorem \ref{main_thm}). This convergence rate results are analogous to the results for gradient descent with noise \cite{ge2015escaping, jin2017escape}. The fundamental difference is that while gradient descent escapes saddles at an exponential rate with a base that depends on the problem's condition number, NCN escapes saddles at an exponential rate with a base of 3/2 for all non-degenerate saddles (Section \ref{sec_illustrative_example}). Section \ref{sec_numerical} considers the problem of matrix factorization to support theoretical conclusions. 

\subsection{Related work} Gradient descent for nonconvex functions converges to an epsilon neighborhood of a critical point, which could be a saddle or a local minimum, in $O(1/\varepsilon^2)$ iterations \cite{nesterov2013introductory}. Escaping saddle points is therefore a fundamental problem for which several alternatives have been developed. A line of work in this regard consists in adding noise when entering a neighborhood of the stationary point. The addition of noise ensures that with high probability the iterates will be at a distance sufficiently large from the stable manifold of the saddle, hence reaching the fundamental conclusion that escape from the saddle point can be achieved in $O(\log(1/\varepsilon))$ \cite{ge2015escaping,jin2017escape} iterations. Noisy gradient descent therefore converges to an epsilon local minimum in $O(1/\varepsilon^2)$ iterations, matching the rate of convergence of gradient descent to stationary points. Under assumptions of nondegeneracy the iterations needed to converge to a local minimum is $O(\log(1/\varepsilon))$. Although the rate of convergence and the rate of escape from saddles match the corresponding rates for NCN, NCN escapes saddles with an exponential base 3/2 but gradient descent escapes saddles with an exponential rate dependent on the condition number. This difference is very significant in practice (Sections \ref{sec_illustrative_example} and \ref{sec_numerical}).

A second approach to ensure that the stationary point attained by the local descent method is a local minimum utilizes second order information to guarantee that the stationary point is a local minimum. These include cubic regularization \cite{griewank1981modification, nesterov2006cubic, cartis2011adaptive, cartis2011adaptive2, agarwal2016finding} and trust region algorithms \cite{conn2000trust, dauphin2014identifying,  curtis2014trust}, as well as approaches where the descent is computed only along the direction corresponding to the negative eigenvalues of the Hessian \cite{curtis2017exploiting}. When using a cubic regularization of a second order approximation of the objective the number of iterations needed to converge to an epsilon local minimum can be further shown to be of order $O(1/\varepsilon^{1.5})$ \cite{nesterov2006cubic}. Solving this cubic regularization is in itself computationally prohibitive. This is addressed with trust region methods that reduce the computational complexity and still converge to a local minimum in $O(1/\varepsilon^{1.5})$ iterations \cite{curtis2014trust}. A related attempt utilizes low-complexity Hessian-based accelerations to achieve convergence in $O(1/\varepsilon^{7/4})$ iterations \cite{agarwal2016finding, carmon2016accelerated}. Although these convergence rates seem to be worse than the rate $O(\log(1/\varepsilon))$ achieved by NCN this is simply a difference in assumptions because we assume here that saddles are nondegenerate. This assumption is absent from \cite{griewank1981modification, nesterov2006cubic, cartis2011adaptive, cartis2011adaptive2, agarwal2016finding, conn2000trust, dauphin2014identifying,  curtis2014trust, curtis2017exploiting}. {If there were degenerate saddles our algorithm would converge to one of them in $O(1/\varepsilon^2)$ iterations as well. It is also worth pointing out that some problems like empirical risk minimization \cite{mei2016landscape} do satisfy the non-degeneracy assumption.}


\section{Nonconvex Newton Method (NCN)}\label{sec_problem_formulation}
Given a multivariate nonconvex function $f: \mathbb{R}^n\to \mathbb{R}$, we would like to solve the following problem 
\begin{equation}\label{eqn_global_minimum}
   \bbx^* := \argmin_{\bbx \in \mathbb{R}^n}\ f(\bbx).
\end{equation}
Finding $\bbx^*$ is NP hard in general, except in some particular cases, e.g., when all local minima are known to be global. We then settle for the simpler problem of finding a local minima $\bbx^\dagger$, which we define as any point where the gradient is null and the Hessian is positive definite 
\begin{equation}\label{eqn_local_minima}
   \left\|\nabla f(\bbx^\dagger) \right\| = 0, \qquad  
   \nabla^2 f(\bbx^\dagger) \succ 0.
\end{equation}
The fundamental difference between (strongly) convex and nonconvex optimization is that any local minimum is global because there is only one point at which $\left\|\nabla f(\bbx) \right\| = 0$ and that point satisfies $\nabla^2 f(\bbx) \succ 0$. Nonconvex functions may have many minima and many other critical points at which $\left\|\nabla f(\bbx^\dagger) \right\| = 0$ but the Hessian is not positive definite. Of particular significance are saddle points, which are defined as those at which the Hessian is indefinite
\begin{equation}\label{eqn_saddles}
   \left\|\nabla f(\bbx^\ddagger) \right\| = 0, \qquad 
   \nabla^2 f(\bbx^\ddagger) \nsucc 0,          \qquad 
   \nabla^2 f(\bbx^\ddagger) \nprec 0.
\end{equation}
Local minima can be found with local descent methods. The most widely used of these is gradient descent which can be proven to approach some $\bbx^\dagger$ with probability one relative to a random initialization under some standard regularity conditions \cite{lee2016gradient,panageas2016gradient}. Convergence guarantees notwithstanding, gradient descent methods can perform poorly around saddle points. Indeed, while escaping saddles is guaranteed in theory, the number of iterations required to do so is large enough that gradient descent can converge to saddles in practical implementations \cite{pascanu2014saddle}. 

Newton's method ameliorates slow convergence of gradient descent by premultiplying gradients with the Hessian inverse. Since the Hessian is positive definite for strongly convex functions, Newton's method provides a descent direction and converges to the minimizer at a quadratic rate {in a neighborhood of the minimum}. The reason for the improvement in the convergence of Newton's method as compared with gradient descent is due to the fact that by premultiplying the descent direction by the inverse of the Hessian we are performing a local change of coordinates by which the level sets of the function become circular. The algorithm proposed here relies in performing an analogous transformation that turns saddles with ``slow'' unstable manifolds as compared to the stable manifold -- this is smaller absolute value of the negative eigenvalues of the Hessian than its positive eigenvalues-- into saddles that have the same absolute values of every eigenvalue. For nonconvex functions the Hessian is not necessarily positive definite and convergence to a minimum is not guaranteed by Newton's method. In fact, all critical points are stable relative to Newton dynamics and the method can converge to a local minimum, a saddle or a local maximum. This shortcoming can be overcome by adopting a modified inverse using the absolute values of the Hessian eigenvalues \cite{nocedal2006numerical}.
%
%
\smallskip\begin{definition}[\bf PT-inverse]\label{def_trucho_inverse}
Let $\bbA\in\reals^{n\times n}$ be a symmetric matrix, $\bbQ\in\reals^{n\times n}$ a basis of orthonormal eigenvectors of $\bbA$, and $\bbLam\in\reals^{n\times n}$ a diagonal matrix of corresponding eigenvalues. %
We say that $\left| \bbLam\right|_m \in \reals^{n\times n}$ is the Positive definite Truncated (PT)-eigenvalue matrix of $\bbA$ with parameter $m$ if 
\begin{equation}
   \left(\left| \bbLam\right|_m\right)_{ii} 
      = \left\{ \begin{array}{l l} 
              \left|\bbLam_{ii}\right| \quad &\mbox{if}\quad \left| \bbLam_{ii}\right| \geq m \\
              m \quad &\mbox{otherwise}.
        \end{array} \right. 
\end{equation}
The PT-inverse of $\bbA$ with parameter $m$ is the matrix $\left| \bbA \right|_m^{-1} = \bbQ \left| \bbLam \right|_m^{-1} \bbQ^{\top}.$
\end{definition}

%

Given the decomposition $\bbA=\bbQ\bbLam\bbQ^{\top}$, the inverse, when it exists, can be written as $\bbA^{-1}=\bbQ\bbLam^{-1}\bbQ^{\top}$. The PT inverse $\left| \bbA \right|_m^{-1} = \bbQ \left| \bbLam \right|_m^{-1} \bbQ^{\top}$ flips the signs of the negative eigenvalues and truncates small eigenvalues by replacing $m$ for any eigenvalue with absolute value smaller than $m$. Both of these properties are necessary to obtain a convergent Newton method for nonconvex functions. We use the PT-inverse of the Hessian to define the NCN method. To do so, consider iterates $\bbx_k$, a step size $\eta_k >0$, and use the shorthand $\bbH_m(\bbx_k)^{-1}=\left|\nabla^2 f(\bbx_k)\right|_m^{-1}$ to represent the PT-inverse of the Hessian evaluated at the $\bbx_k$ iterate. The NCN method is defined by the recursion
\begin{equation}\label{eqn_update}
   \bbx_{k+1}\ =\ \bbx_k - \eta_k \bbH_m(x_k)^{-1}\nabla f(\bbx_k)
             \ =\ \bbx_k - \eta_k \left|\nabla^2 f(\bbx_k)\right|_m^{-1}\nabla f(\bbx_k).
\end{equation}
The step size $\eta_k$ is chosen with a backtracking line search as is customary in regular Newton's method; see, e.g., \cite[Section 9.5.2]{boyd2004convex}. This yields a step routine that is summarized in Algorithm \ref{alg_ncn}. In Step 3 we update the iterate $\bbx_k$ using the PT-inverse Hessian $\bbH_m(\bbx_k)^{-1}$ computed in Step 2 and initial stepsize $\eta_k=1$. The updated variable $\bbx_{k+1}$ is checked against the decrement condition with parameter $\alpha\in(0,0.5)$ in Step 4. If the condition is not met, we decrease the stepsize $\eta_k$ by backtracking it with the constant $\beta<1$ as in Step 5. We update the iterate $\bbx_k$ with the new stepsize as in Step 6 and repeat the process until the decrement condition is satisfied.

%
{\linespread{1.4}
\begin{algorithm}[t] \begin{algorithmic}[1]
    \STATE \textbf{function} $[\bbx_{k+1}]=$ NCNstep$(\bbx_k,\alpha,\beta)$
    
\STATE Compute PT-inverse: $\bbH_m(\bbx_k)^{-1}=\left|\nabla^2 f(\bbx_k)\right|_m^{-1}$.
\STATE Set step size to $\eta_k=1$. Update argument to $\bbx_{k+1} = \bbx_k - \eta_k \bbH_m(\bbx_k)^{-1}\nabla f(\bbx_k)$.
   \WHILE[Backtracking]{$f(\bbx_{k+1})>f(\bbx_{k})-\alpha\eta_k\nabla f(\bbx_k)^\top\bbH_m(\bbx_k)^{-1}\nabla f(\bbx_k)$} 
   \STATE Reduce step size to $\eta_k = \beta \eta_k$.
   \STATE Update argument to $\bbx_{k+1} = \bbx_k - \eta_k \bbH_m(\bbx_k)^{-1}\nabla f(\bbx_k)$.
   \ENDWHILE\ \{return $\bbx_{k+1}$\}
\end{algorithmic}
\caption{Non-Convex Newton Step}\label{alg_ncn} \end{algorithm}}

%
Since the PT-inverse is defined to guarantee that $\bbH_m(x_k)^{-1}\nabla f(\bbx_k)$ is a proper descent direction, it is unsurprising that NCN converges to a local minimum. The expectation is, however, that it will do so at a faster rate because of the Newton-like correction that is implied by \eqref{eqn_update}. Intuitively, the Hessian inverse in convex functions implements a change of coordinates that renders level sets approximately spherical around the current iterate $\bbx_k$. The Hessian PT-inverse in nonconvex functions implements an analogous change of coordinates that renders level sets in the neighborhood of a saddle point close to a symmetric hyperboloid. This regularization of level sets is expected to improve convergence, something that has been observed empirically, \cite{pascanu2014saddle}.


\subsection{Convergence of NCN to local minima}\label{sec_convergence_preview}

Convergence results are derived with customary assumptions on Lipschitz continuity of the gradient and Hessian, boundedness of the norm of the local minima, and non-degeneracy of critical points:

%
\begin{assumption}\label{assumption_lipschitz}
The function $f(\bbx)$ is twice continuously differentiable. The gradient and Hessian of $f(\bbx)$ are Lipchitz continuous, i.e., there exits constants $M, L>0$ such that for any $\bbx, \bby \in \mathbb{R}^n$
  \begin{equation}
\left\|\nabla f(\bbx) - \nabla f(\bby)\right\| \leq M \| \bbx-\bby\|, \quad \mbox{and} \quad 
  \left\|\nabla^2 f(\bbx) - \nabla^2 f(\bby)\right\| \leq L \| \bbx-\bby\|.
  \end{equation}  
\end{assumption}

\begin{assumption}\label{assumption_bounded_opt}
There exists a positive constant $B$ such that the norm of local minima satisfies $\|\bbx^\dagger\|\leq B$ for all $\bbx^\dagger$ satisfying \eqref{eqn_local_minima}. In particular, this is true of the global minimum $\bbx^*$ in \eqref{eqn_global_minimum}.
\end{assumption}

{
\begin{assumption}\label{assumption_morse} Local minima and saddles are non-degenerate. I.e., there exists a constant $\xi>0$ such that $\min_{i=1\ldots n} \lambda_i\left(\nabla^2 f(\bbx^\dagger)\right)>\xi$ for all local minima $\bbx^\dagger$ and $\min_{i=1\ldots n} \left|\lambda_i\left(\nabla^2 f(\bbx^\ddagger)\right)\right| >\xi$ for all saddle pioints $\bbx^\ddagger$ defined in \eqref{eqn_saddles}. The notation $\lambda_i(\nabla^2f(\bbx))$ refers to the $i$-th eigenvalue of the Hessian of $f$ at the point $\bbx$.
\end{assumption}
}
%
The main feature of the update in \eqref{eqn_update} is that it exploits curvature information to accelerate the rate for escaping from saddle points relative to gradient descent. In particular, the iterates of NCN escape from a local neighborhood of saddle points exponentially fast at a rate which is independent of the problem's condition number. {This neighborhood is defined as $\left\|\nabla f(\bbx)\right\|<\delta/2$ where
\begin{equation}\label{def_delta}
  \delta = \min\left\{ \frac{m^2\left(1-2\alpha\right)}{L}, \frac{m^2}{5L}\right\}.
\end{equation}
Throughout the paper, we make the assumption that the accuracy with which we want to solve the problem is $\varepsilon$ satisfies $\varepsilon<\delta$. 
}
To state this result formally, let $\bbx^{\ddagger}$ be a saddle of interest and denote $\bbQ_{-}$ and $\bbQ_{+}$ as the orthogonal subspaces associated with the negative and positive eigenvalues of $\nabla^2 f(\bbx^\ddagger)$. For a point $\bbx\neq\bbx^{\ddagger}$ we define the gradient projections on these subspaces as
\begin{equation}\label{def_grad_neg_and_pos}
   \nabla f_{-}(\bbx) := \bbQ_{-}^\top \nabla f(\bbx) \quad \mbox{and} \quad 
   \nabla f_{+}(\bbx) := \bbQ_{+}^\top \nabla f(\bbx).
\end{equation}
These projections have different behaviors in the neighborhood of a saddle point. The projection on the positive subspace $\nabla f_{+}(\bbx) $ enjoys an approximately quadratic convergent phase as in Newton's method (Theorem \ref{theo_rate_saddles}). This is as would be expected because the positive portion of the Hessian is not affected by the PT-inverse. The negative portion $ \nabla f_{-}(\bbx)$ can be shown to present an exponential divergence from the saddle point with a rate independent of the problem conditioning. These results provide a bound in the number of steps required to escape the neighborhood of the saddle point that we state next. 
\begin{theorem}\label{thm_steps_escape}
  Let $f(\bbx)$ be a function satisfying Assumptions \ref{assumption_lipschitz} and \ref{assumption_morse}, $\varepsilon>0$ be the desired accuracy of the solution provided by Algorithm \ref{alg_nonconvex_newton} and $\alpha \in (0,1)$ be one of its inputs. If {$m<\xi/2$ and}
  \begin{equation}
    \left\|\nabla f_{-}(\bbx_0) \right\| \geq \max\{({5L}/{2m^2})\left\|\nabla f(\bbx_0) \right\|^2,\varepsilon\}
  \end{equation}
    and $ \left\|\nabla f(\bbx_0) \right\| \leq \delta/2$, we have that $\left\|\nabla f(\bbx_{K_1}) \right\| \geq \delta/2$, with 
\begin{equation}\label{eqn_number_steps_escape}
K_1 \leq 1+\log_{3/2}\left( \frac{\delta}{2\varepsilon}\right).
\end{equation}
\end{theorem}
\begin{proof}
  {This theorem is a particular case of the Theorem \ref{theo_rate_saddles} in
    Section \ref{sec_analysis}.}
  \end{proof}
The result in Theorem \ref{thm_steps_escape} establishes an upper bound for the number of iterations needed to escape the saddle point which is of the order $O(log(1/\varepsilon))$ as long as the iterate $\bbx_0$ satisfies $\left\|\nabla f_{-}(\bbx_0) \right\| \geq ((5L)/(2m^2))\left\|\nabla f(\bbx_0) \right\|^2$ and $ \left\|\nabla f_{-}(\bbx_0) \right\| \geq \varepsilon$. However, the fundamental result is that the rate at which the iterates escape the neighborhood of the saddle point is a constant $3/2$ independent of the constants of the specific problem. To establish convergence to a local minimum we will prove four additional results:

\noindent (i) In Proposition \ref{prop_sublinear} we state that the convergence of the algorithm to a neighborhood of the critical points such that $\left\|\nabla f(\bbx_k)\right\| <\delta/2$ is achieved in a constant number of iterations bounded by
\begin{equation}\label{eqn_number_steps_sublinear}
  K_2 = \frac{4M^2\left(f(\bbx_0) -f(\bbx^*)\right)}{\alpha \beta m \delta^2}.
\end{equation}
(ii) We show in Proposition \ref{prop_revisit_saddle} that the number of times that the iterates re-visit the same neighborhood of a saddle point is upper bounded by 
\begin{equation}
  T<\frac{2}{\alpha\beta}\frac{M^3}{m^3}+\alpha\beta,
  \end{equation}
and that (iii) once in such neighborhood of a local minimum, the algorithm achieves $\varepsilon$ accuracy in a number of iterations bounded by (Corollary \ref{coro_local_minima}) 
\begin{equation}\label{eqn_number_steps_local_minima}
  K_3=\log_2\left(\log_2\left(\frac{2m^2}{5L\varepsilon}\right)\right).
\end{equation}
(iv) For the case that the iterate is within the neighborhood of a saddle point, but the conditions required by Theorem \ref{thm_steps_escape} are not satisfied, we show that by adding noise to the iterate we can ensure that said conditions are met with probability $(1-p)^{1/S}$ after a number of iterations of order $O(\log(\log(1/\varepsilon))+\log(S/p))$, where $S$ is the {maximum number of saddles that the algorithm visits. To converge to the minimum with probability $1-p$ we need to escape each one of them with probability $(1-p)^{1/S}$.} In particular, we show that if we add a bounded version of the Gaussian noise  $ \ccalN(0, 2\varepsilon/m)$ to each component of the decision variable, with probability $q$ the perturbed variable will be in the region that conditions required by Theorem \ref{thm_steps_escape} are satisfied. We further show that the probability $q$ is lower bounded by 
\begin{equation}\label{eqn_expression_q}
q>2\left(1-\Phi(1) \right)\frac{\gamma\left(\frac{n}{2},\frac{n}{2}\right)}{\Gamma\left(\frac{n}{2}\right)},
\end{equation}
  where $\Phi(1)$ is the integral of the Gaussian distribution with integration boundaries $-\infty$ and $1$, $\Gamma(\cdot)$ is the gamma function, $\gamma(\cdot,\cdot)$ is the lower incomplete gamma function and $n$ is the dimension of the iterate $\bbx$. {In practice however, we cannot verify if the conditions required in Theorem \ref{thm_steps_escape} are satisfied by the perturbed iterate. To solve this issue, we show that the probability $q$ defined in \eqref{eqn_expression_q} is in fact the probability of falling in the region $\left\|\nabla f_{-}(\bbx_0) \right\| \geq \max\{({5L}/{2m^2})\left\|\nabla f(\bbx_0) \right\|^2,\varepsilon/2\}$. In this region, the iterates present the same behavior than in Theorem \ref{thm_steps_escape}, and therefore the number of steps required to escape the $\varepsilon$ neighborhood is $1+\log_{3/2}\left(2\varepsilon/\varepsilon\right)<2$. Hence, after the perturbation, we run the algorithm for $2$ iterations, i.e., the maximum number of required iterations to escape the $\varepsilon$ neighborhood of the saddle assuming that the perturbed iterate is in the preferable region.} If the perturbed variable were not in the region of interest, the iterates may not escape the saddle and another round of perturbation is needed. Using this scheme we show that the iterates will be within the desired region of the saddle with probability $(1-p)^{1/S}$ after at most 
  \begin{equation}\label{eqn_number_steps_rand}
  \begin{split}
 K_4 \leq \left(1+\frac{\log (S/p)}{ \log \left(1/(1-q)\right)}\right)\left[ \log_2 \left(  \log_2\left(\frac{5L}{2m^2\varepsilon}\right)\right)+\log_{3/2}(2)+1\right] .
    \end{split}
  \end{equation}
  The fact that we may need to visit each saddle $T$ times does not contribute to the increase in the previous probability as we show in Proposition \ref{prop_revisit_saddle}, since in only one of theese visits we reach the $\varepsilon$ neighborhood, where noise is added. By combing the previous bounds we establish a total complexity of order $O(\log(1/p)+\log(1/\varepsilon))$ for NCN to converge to an $\varepsilon$ neighborhood of a local minima of $f$ with probability $1-p$. We formalize this result next.

{\linespread{1.3}
\begin{algorithm}[t] \begin{algorithmic}[1]
    \STATE \textbf{Input:} $\bbx_k = \bbx_0$, accuracy $\varepsilon>0$ and parameters $\alpha\in(0,1/2)$, $\beta\in (0,1)$ 
    
\WHILE {$\left\|\nabla f(\bbx_k) \right\|>\varepsilon$ or $\min \lambda(\nabla^2 f(\bbx_k))<0$} 
\STATE Compute the updated variable $\bbx_{k+1}=$ NCNstep$(\bbx_k,\alpha,\beta)$

   \IF {$\left\|\nabla f(\bbx_{k+1}) \right\| \leq \varepsilon$ and $\min \lambda(\nabla^2 f(\bbx_{k+1}))<0$}
   \STATE{  $\tilde{\bbx} = \bbx_{k+1} +X$} with {$X_i \sim \ccalN(0,2\varepsilon/m)$}
   \WHILE{$\left\|\nabla f(\tilde{\bbx})\right\|> (2\sqrt{n}M/m+1)\varepsilon$}
   \STATE{  $\tilde{\bbx} = \bbx_{k+1} +X$} with {$X_i \sim \ccalN(0,2\varepsilon/m)$}
   \ENDWHILE
   \STATE{$\bbx_{k+1}=\tilde{\bbx}$}
   	\IF {$\left\| \nabla f(\bbx_{k+1})\right\| \leq \varepsilon$} 
	\STATE Update the iterate $\bbx_{k+1}=$ NCNstep$(\bbx_k,\alpha,\beta)$ for $\lceil\log_{3/2}2\rceil = 2$ times
	\ENDIF
   \ENDIF
\ENDWHILE
\end{algorithmic}
\caption{Non-Convex Newton Method}\label{alg_nonconvex_newton} \end{algorithm}}
%
%
%
\begin{theorem}\label{main_thm}
  Let $f(\bbx)$ be a function, satisfying Assumptions \ref{assumption_lipschitz}-- \ref{assumption_morse}. Let $\varepsilon>0$ be theaccuracy of the solution and $\alpha\in(0,1/2)$, $\beta\in (0,1)$ the remaining inputs of Algorithm \ref{alg_nonconvex_newton}. Then if {$m<\xi/2$}, with probability $1-p$ and for any $\varepsilon$ satisfying
  \begin{equation}
  \varepsilon < \min \left\{\frac{1}{8Ln},\frac{\delta m}{4M\sqrt{n}+m},\frac{L}{5\left(2M\sqrt{n}+m\right)^2},\sqrt{\frac{\alpha\beta\frac{m^3}{M^3}( \delta/2)^2}{1+\left(2\sqrt{n}\frac{M}{m}+1\right)^2}}\right\},
    \end{equation}
Algorithm \ref{alg_nonconvex_newton} outputs $\bbx_{sol}$ satisfying $\left\|\nabla f(\bbx_{sol}) \right\| < \varepsilon$ and $\nabla^2 f(\bbx_{sol}) \succ 0$ in at most 
 %
 %
 \begin{equation}\label{eqn_final_complexity}
   \begin{split}
     K=STK_1+(ST+1)K_2 +K_3+SK_4
\end{split}
 \end{equation}
%
 iterations, where $K_1$, $K_2$, $K_3$, $K_4$ and $T$ are the bounds are defined in \eqref{eqn_number_steps_escape}--\eqref{eqn_number_steps_rand} {and $S$ -- the maximum number of saddles that the algorithm approach -- is bounded by
 \begin{equation}
S \leq \frac{m^2}{M\alpha\beta \left(\delta/2\right)^2}\left(f(\bbx_0)-f(\bbx^*)\right).
   \end{equation}
}
\end{theorem}
The result in Theorem \ref{main_thm} states that the proposed NCN method outputs a solution with a desired accuracy and with a desired probability in a number of steps that is bounded by $K$. The final bound follows from the fact that it may be required to visit every saddle point $T$ times before reaching a local minima. Hence we may have  $TS+1$ approaches to neighborhoods of the critical points, each one of them taking $K_2$ iterations. The latter corresponds to the second term in \eqref{eqn_final_complexity}. The first term correspond to the need of escaping all $S$ saddles, $T$ times each if we are if we are in the good region. Thus taking a total of $STK_1$ steps. If we are not in the desired region, it takes $K_4$ steps to reach said region and we may have to do this in all $S$ saddles, hence the last term in \eqref{eqn_final_complexity}.Finally, term $K_3$ correspond to the quadratic convergence to the local minimum. The dominant terms in \eqref{eqn_final_complexity} are $K_1$, which has a $O(\log(1/\varepsilon))$ dependence, and $K_4$, which depends on the probability as $O(\log(1/p))$. 

Before proving the result of the previous theorem we describe the details of Algorithm \ref{alg_nonconvex_newton}. Its main core is the NCN step described in \eqref{eqn_update} (Step 3) that is performed as long as the iterates are not in a neighborhood of a local minimum satisfying $\left\| \nabla f(\bbx_k)\right\|<\varepsilon$. Steps 4--12 are introduced to add Gaussian noise to satisfy the hypothesis of Theorem \ref{thm_steps_escape}. If the iterate is in a neighborhood of a saddle point such that $\left\| \nabla f(\bbx_k)\right\| <\varepsilon$ (Step 4), noise from a Gaussian distribution is added (Step 5). The draw is repeated as long as $\left\| \nabla f(\bbx_k)\right\| >(2\sqrt{n}M/m+1)\varepsilon$ (Steps 6--8). This is done to ensure that the iterates remain close to the saddle point. Once this condition is satisfied we perform the NCN step twice if the iterate is still in the neighborhood $\left\| \nabla f(\bbx_k)\right\| < \varepsilon$ (steps 10--12). In cases where $\left\|\nabla f({\bbx_k}) \right\|\geq (5L/2m^2)\left\|\nabla f(\bbx_k)\right\|^2$, two steps of NCN are enough to escape the $\varepsilon$ neighborhood of the saddle point and therefore to satisfy the hypothesis of Theorem \ref{thm_steps_escape} (c.f. Proposition \ref{prop_probabilistic}).

\subsection{An Illustrative Example}\label{sec_illustrative_example}
%
%
\begin{figure}
    \centering
        \includegraphics[width=\textwidth, height=0.62\textwidth]{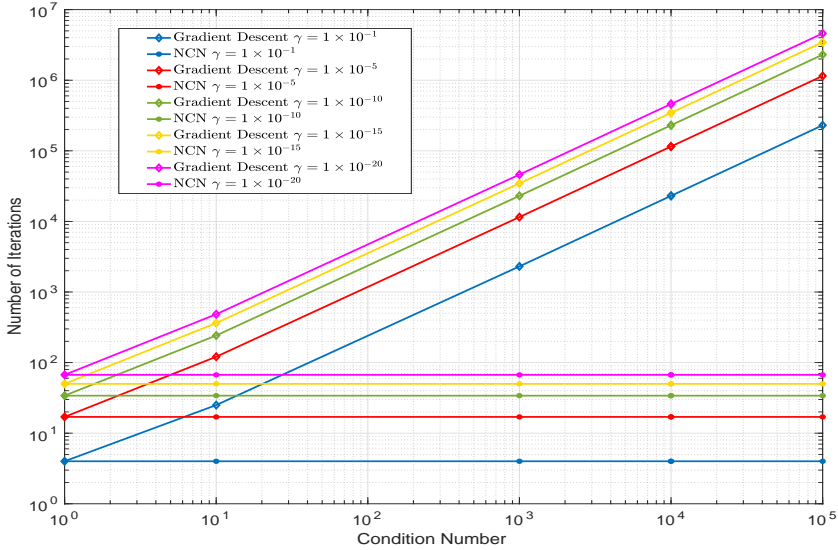}
        \caption{Number of iterations to escape the neighborhood $[-1,1]^2$ for quadratic problems with different condition numbers $\lam^{-1}$ and different initial distance to the stable manifold $\gamma$. NCN takes the same number of steps to escape the neighborhood of the saddle point regardless of the condition number of the problem. In the case of gradient descent this dependence is linear.} 
        \label{example_escape}
\end{figure}
In understanding escape from a saddle point it is important to distinguish challenges associated with the saddle's condition number and challenges associated with starting from an initial point that is close to the stable manifold. To illustrate these two different challenges we consider a family $f_{\lam}$ of nonconvex functions parametrized by a coefficient $\lambda\in(0,1]$ and defined as
\begin{equation}
   f_{\lam}(\bbx)\ =\ \frac{1}{2}x_1^2- \frac{\lambda}{2}x_2^2,
\end{equation}
As $\lam$ decreases from 1 to 0 the saddle becomes flatter, its condition number $\lam^{-1}$ growing from $1$ to $+\infty$. Gradient descent iterates using unit stepsize for the function $f_{\lam}$ result in zeroing of the first coordinate in just one step because $x_1^+ =  x_1-\eta \nabla_{x_1} f_\lam(\bbx) = x_1 - x_1 = 0$. The second coordinate evolves according to the recursion
\begin{equation}
   x_2^+\ =\  x_2-\eta \nabla_{x_2} f_\lam(\bbx)\ =\ (1+{\lambda})x_2
\end{equation}
Likewise, NCN with unit step size results in zeroing of the first coordinate because $x_1^+ = x_1-\eta [\bbH^{-1}\nabla f_\lam(\bbx)]_1 = 0$. The second coordinate, however, evolves as
\begin{equation}
   x_2^+\ =\ x_2- \left[\bbH^{-1}\nabla f_\lam(\bbx)\right]_2\ =\ 2x_2.
\end{equation}
Both expressions imply exponential escape from the saddle point. In the case of gradient descent the base of escape is $(1+\lam)$ but in the case of NCN the base of escape is 2 independently of the value of $\lam$.

The consequences of this difference are illustrated in Figure \ref{example_escape} which depicts the number of iterations it takes to reach the border of the $L_\infty$ unit ball as a function of the condition number $\lam^{-1}$. Different curves depict different initial conditions in terms of their distance $\gamma$ to the stable manifold $x_2=0$, which is simply the initial value for the $x_2$ coordinate. The escape from the saddle for NCN occurs in $\ln\gamma / \ln 2$ iterations, a number that is independent of the condition number. The time it takes for gradient descent to escape a saddle $\ln\gamma / \ln (1+\lam) \approx \lam^{-1}\ln\gamma$. This number is roughly proportional to the saddle's condition number. An important observation that follows from Figure \ref{example_escape} is that the condition number of the saddle is a more challenging problem that the initial distance to the stable manifold. If the saddle is well conditioned, escape from the saddle with gradient descent takes a few iterations even when the initial condition is very close to the stable manifold. E.g., $\ln\gamma / \ln 2 < 67$ when $\gamma=10^{-20}$. If the saddle is not well conditioned escape from the saddle with gradient descent takes a very large number of iterations even if the initial condition is far from the stable manifold. E.g.,  $\ln\gamma / \ln (1+\lam)>10^5$ when $\gamma=0.1$ and $\lam=10^{-5}$. This is because the number of iterations to escape a saddle is approximately $\lam^{-1}\ln\gamma$. This number grows logarithmically with $\gamma$ but linearly with the condition number $\lam^{-1}$. In the case of NCN, escape takes always $\ln\gamma / \ln 2 $ iterations independently of $\lam$. Theorem \ref{thm_steps_escape} provides a qualitatively analogous statement for any saddle.

\section{Convergence Analysis}\label{sec_analysis}
%
To study the convergence of the proposed NCN method we divide the results into two parts. In the first part, we study the performance of NCN in a neighborhood of critical points. We first define this region and then characterize the number of required iterations to escape it in the case that the critical point is a saddle or the number of required iterations for convergence in case that the critical point is a minimum. In the second part, we study the behavior of NCN when the iterates are not close to a critical point and derive an upper bound for the number of iterations required to reach one.

To analyze the local behavior of NCN we characterize the region in which the step size chosen by backtracking line search is $\eta=1$, as in standard Newton's method for convex optimization. We formally introduce this region in the following lemma.

\begin{lemma}\label{lemma_step_size_one}
    Let $f$ be a function satisfying Assumptions \ref{assumption_lipschitz} and \ref{assumption_morse}, and $\alpha \in (0,1/2)$ be the input parameter of Algorithm \ref{alg_nonconvex_newton}. Then if $\|\nabla f(\bbx)\|\leq ({m^2}/{L})\left(1-2\alpha\right)$, a backtracking algorithm admits the step size $\eta=1$.
\end{lemma}
\begin{proof}
See Appendix \ref{sec_step_size_one}.
\end{proof}
%
The result in Lemma \ref{lemma_step_size_one} characterizes the neighborhood in which the step size of NCN chosen by backtracking is $\eta=1$. In the following theorem, we study the behavior of NCN in this neighborhood. Before introducing the result, recall the definitions of the gradient projected over the subspace of the eigenvectors associated with the negative and positive eigenvalues in \eqref{def_grad_neg_and_pos} which are denoted by $\nabla f_{-} (\bbx)$ and $\nabla f_{+} (\bbx)$, respectively. We attempt to show that the norm $\|\nabla f_{-}(\bbx_k)\|$ is almost doubled per iteration in this local neighborhood, and the norm $\|\nabla f_{+} (\bbx)\|$ converges to zero quadratically. 
%
\begin{theorem}\label{theo_rate_saddles}
  Let Assumptions \ref{assumption_lipschitz}--\ref{assumption_morse} and $\left \|\nabla f(\bbx_k)\right\| < m^2(1-2\alpha)/L$ hold, \blue{$m<\xi/2$}. Then $\nabla f_{-}(\bbx)$ and $\nabla f_{+}(\bbx)$ defined in \eqref{def_grad_neg_and_pos} satisfy
  \begin{equation}\label{eqn_rate_saddles_pos}
\| \nabla f_{+} (\bbx_{k+1})\| \leq \frac{5L}{2m^2} \left\|\nabla f( \bbx_k) \right\|^2,
  \end{equation}
  and
  \begin{equation}\label{eqn_rate_saddles_neg}
\left\|    \nabla f_{-}(\bbx_{k+1}) \right\| \geq 2\left\|\nabla f_{-}(\bbx_k)\right\| -\frac{5L}{2m^2} \left\|\nabla f( \bbx_k) \right\|^2.
  \end{equation}
\end{theorem}
\begin{proof}
It follows from Lemma \ref{lemma_step_size_one} that the backtracking algorithm admits a step $\eta=1$. Hence we can write $\nabla f(\bbx_{k+1})$ as
  \begin{equation}
\nabla f(\bbx_{k+1}) = \nabla f(\bbx_k) + \int_0^1\nabla^2f(\bbx_k +\theta \Delta \bbx)\Delta \bbx d\theta. 
  \end{equation}
We next show that in the region $\left\| \nabla f(\bbx_k)\right\| < m^2(1-2\alpha)/L$ it holds that $m<\min_{i=1\ldots n} \left|\lambda_i\left(\nabla^2 f(\bbx)\right)\right|$. Let us consider the region $\left\| \bbx_k-\bbx_c\right\|< m(1-2\alpha)/L$. In this region, by virtue of Lemma \ref{lemma_bound_eigenvalues} we have that $\min_{i=1\ldots n} \left|\lambda_i\left(\nabla^2 f(\bbx)\right)\right|>m$ and that in the boundary $\left\| \nabla f(\bbx_k) \right\| \geq m^2(1-2\alpha)/L$. Since at the critical point $\left\| \nabla f(\bbx_c)\right\| = 0$, by continuity of the norm we have that the region $\left\| \nabla f(\bbx_k)\right\| < m^2(1-2\alpha)/L$ is contained in the region in which $\min_{i=1\ldots n} \left|\lambda_i\left(\nabla^2 f(\bbx)\right)\right|>m$. Thus,  we have that $\bbH_m(\bbx_k) = \bbQ(\bbx_k)\left|\Lambda(\bbx_k) \right| \bbQ(\bbx_k)^\top$. Using the fact that $\nabla f(\bbx_k) = \bbH_m(\bbx_k)\bbH_m(\bbx_k)^{-1} \nabla f(\bbx_k) = -\bbH_m(\bbx_k)\Delta \bbx$ the previous expression can be written as
\begin{equation}\label{eqn_gradient_expansion}
  \begin{split}
    \nabla f(\bbx_{k+1}) &=2\nabla f(\bbx_k) + \int_0^1\left(\nabla^2f(\bbx_k+\theta\Delta \bbx)+\bbH_m(\bbx_k)\right)\Delta \bbx d\theta \\
    &=2\nabla f(\bbx_k)  +\int_0^1\left(\nabla^2f(\bbx_k+\theta\Delta \bbx)-\nabla^2 f(\bbx_k)\right)\Delta \bbx d\theta\\
    &+\int_0^1\left(\nabla^2f(\bbx_k)-\nabla^2 f(\bbx_c)\right)\Delta \bbx d\theta  +\int_0^1\left(\nabla^2f(\bbx_c)+\bbH_m(\bbx_c)\right)\Delta \bbx d\theta \\
&     +\int_0^1\left(\bbH_m(\bbx_k)-\bbH_m(\bbx_c)\right)\Delta \bbx d\theta.
\end{split}
  \end{equation}
The latter equality follows by adding and substracting $\nabla^2 f(\bbx_c)\Delta \bbx$, $\nabla^2 f(\bbx_k)\Delta \bbx$ and $\bbH_m(\bbx_c)\Delta \bbx$. Multiply both sides of \eqref{eqn_gradient_expansion} by $\bbQ_{-}^\top$ the matrix of the eigenvectors corresponding to negative eigenvalues of the Hessian at $\bbx_c$. Since $\bbQ_{-}$ is a matrix whose columns are eigenvectors its norm is bounded by one. Combine this fact with the Lipschitz continuity of the Hessian (c.f. Assumption \ref{assumption_lipschitz}) to write
\begin{equation}\label{eqn_gradient_first_bound}
  \begin{split}
    \left\|\bbQ_{-}^\top \left(\nabla ^2 f(\bbx_k+\theta\Delta \bbx)-\nabla ^2 f(\bbx_k)\right)\Delta\bbx \right\| 
    &\leq L\theta \left\|\Delta\bbx\right\|^2,
\end{split}
  \end{equation}
Likewise, we can upper bound the second and fourth integrands in \eqref{eqn_gradient_expansion} by
\begin{equation}\label{eqn_gradient_second_bound}
\left\|\bbQ_{-}^\top \left(\nabla ^2 f(\bbx_k)-\nabla ^2 f(\bbx_c)\right)\Delta\bbx \right\| \leq L \left\| \bbx_k - \bbx_c\right\| \left\|\Delta \bbx \right\|, 
\end{equation}
\begin{equation}\label{eqn_gradient_third_bound}
\left\|\bbQ_{-}^\top \left(\bbH_m (\bbx_k)-\bbH_m(\bbx_c)\right)\Delta\bbx \right\| \leq L \left\| \bbx_k - \bbx_c\right\| \left\|\Delta \bbx \right\|.
\end{equation}
We next show that the third integrand in \eqref{eqn_gradient_expansion} when multiplied by $\bbQ_{-}^\top$ becomes zero. Let us write the product $\bbQ_{-}^\top \left(\nabla^2 f(\bbx_c) + \bbH_m(\bbx_c)\right)$ as 
\begin{equation}
\bbQ_{-}^\top \left(\nabla^2 f(\bbx_c) + \bbH_m(\bbx_c)\right) = \bbQ_{-}^\top\bbQ\left(\mathbf{\Lambda}(\bbx_c)+\left| \mathbf{\Lambda}(\bbx_c)\right|\right) \bbQ^\top.
  \end{equation}
Let $d$ be the number of negative eigenvalues. Then $\bbQ_{-}^\top \bbQ = \left[\bbI_d, \mathbf{0}_{d\times n-d}\right]$. Moreover $\mathbf{\Lambda}(\bbx_c)+\left| \mathbf{\Lambda}(\bbx_c)\right|$ is diagonal with the first $d$ elements being zero and the remaining $n-d$ being $2$. Which shows that  $\bbQ_{-}^\top \bbQ \left(\mathbf{\Lambda}(\bbx_c)+\left| \mathbf{\Lambda}(\bbx_c)\right|\right)=\mathbf{0}$. With this result and the bounds on \eqref{eqn_gradient_first_bound},\eqref{eqn_gradient_second_bound},\eqref{eqn_gradient_third_bound} we can lower bound \eqref{eqn_gradient_expansion} by
\begin{equation}\label{eqn_gradient_neg_bound}
  \begin{split}
    \left\|    \nabla f_{-}(\bbx_{k+1}) \right\|&\geq 2\left\|\nabla f_{-}(\bbx_k)\right\| -L\left\|\Delta\bbx\right\|^2 \int_0^1\theta d\theta -2L\left\| \bbx_k - \bbx_c\right\| \left\|\Delta \bbx \right\|.
\end{split}
  \end{equation}
Finally, using the fact that $\left\| \nabla f(\bbx_k)\right\| \geq m\left\|\bbx_k-\bbx_c\right\|$ (c.f. Lemma \ref{lemma_bound_eigenvalues}) and that $\left\|\Delta \bbx\right\| \leq \frac{1}{m} \left\|\nabla f(\bbx_k) \right\|$, 
the previous bound reduces to 
\begin{equation}\label{eqn_gradient_neg_final}
  \begin{split}
    \left\|    \nabla f_{-}(\bbx_{k+1}) \right\|&\geq 2\left\|\nabla f_{-}(\bbx_k)\right\| -\frac{5L}{2m^2} \left\|\nabla f( \bbx_k) \right\|^2.
\end{split}
  \end{equation}
The proof for the projection over the positive subspace is analogous.
  \end{proof}
%
The first result in Theorem \ref{theo_rate_saddles} shows that the norm $\| \nabla f_{+} (\bbx_{k+1})\| $ approaches zero at a quadratic rate if most of the energy of the gradient norm $\| \nabla f (\bbx_{k})\|$ belongs to the term $\| \nabla f_{+} (\bbx_{k})\|$. In particular, when we are in a local neighborhood of a local minimizer and all the eigenvalues are positive we have $\| \nabla f (\bbx_{k})\|= \|\nabla f_{+} (\bbx_{k})\|$ and therefore the sequence of iterates converges quadratically to the local minimum. Indeed, in a neighborhood of a local minimum the algorithm proposed here is equivalent to Newton's method. We formalize this result in the following corollary. 
%
\begin{corollary}\label{coro_local_minima}
 Let $f$ be a function satisfying Assumptions \ref{assumption_lipschitz} and \ref{assumption_morse} and let $\bbx_0 \in \mathbb{R}^n$ be in the neighborhood of a local minima such that $\left \|\nabla f(\bbx_0)\right\| < \delta$ where $\delta$ is that in \eqref{def_delta} \blue{and $m<\xi/2$}.
 Then, it holds that $\left\|\nabla f(\bbx_{K_3}) \right\|  <\varepsilon$, where $K_3$ is bounded by
 \begin{equation}
K_3 \leq \log_{2}\left(\log_2\left(\frac{2m^2}{5L\varepsilon}\right)\right).
 \end{equation}
  \end{corollary}
\begin{proof}
 See Appendix \ref{app_local_minima}.
\end{proof}
The second result in Theorem \ref{theo_rate_saddles} shows that the norm $\left\|\nabla f_{-}(\bbx_k)\right\|$ multiplies by factor $2-\zeta$, where $0<\zeta<1$ is a free parameter, after each iteration of NCN if the squared norm $\left\|\nabla f( \bbx_k) \right\|^2$ is negligible relative to $\left\|\nabla f_{-}(\bbx_k)\right\|$. We formally state this result in the following Proposition. 

%
\begin{proposition}\label{prop_steps_escape}
  Let $f$ be a function satisfying Assumptions \ref{assumption_lipschitz} and \ref{assumption_morse}. Further, recall the definition of $\delta$ in \eqref{def_delta} and let $\zeta \in(0,1)$ and $\varepsilon>0$ be the inputs of Algorithm \ref{alg_nonconvex_newton}. {Then, if the conditions $\left\|\nabla f_{-}(\bbx_0) \right\| \geq \max\{(5L)/(2m^2)\left\|\nabla f(\bbx_0) \right\|^2,\varepsilon\}$ and $ \left\|\nabla f(\bbx_0) \right\| \leq \zeta \delta$ hold, \blue{with $m<\xi/2$}, the sequence generated by NCN needs $K_1$ iterations to escape the saddle and satisfy the condition $\left\|\nabla f(\bbx_{K_1}) \right\| \geq \zeta \delta$, where $K_1$ is upper bounded by}
\begin{equation}
K_1 \leq1+ \frac{\log\left( \zeta \delta/\varepsilon\right)}{\log \left(2-\zeta\right)}.
\end{equation}
\end{proposition}
\begin{proof}
  See Appendix \ref{app_steps_escape}.
  \end{proof}
%
The result in Proposition \ref{prop_steps_escape} states that when the norm of the projection of the gradient over the subspace of eigenvectors corresponding to negative eigenvalues $\left\|\nabla f_{-}(\bbx_0) \right\|$ is larger than the required accuracy $\eps$ and the expression $ (5L)/(2m^2)\left\|\nabla f(\bbx_0) \right\|^2$, then the sequence of iterates generated by NCN escapes the saddle exponentially with base $(2-\zeta)$. Since $\zeta$ is a free parameter we set $\zeta=1/2$ as in Theorem \ref{thm_steps_escape}. When the condition $\left\|\nabla f_{-}(\bbx_0) \right\| \geq (5L)/(2m^2)\left\|\nabla f(\bbx_0) \right\|^2$ in Proposition \ref{prop_steps_escape} is not met for all $k$, the sequence generated by NCN reaches the $\varepsilon$ neighborhood after $\log_2 \left( \log_2\left((5L\zeta)/(m^2\varepsilon)\right)\right)$ iterations as shown in the following proposition. 
%
\begin{proposition}\label{prop_quadratic_convergence}
  Let $f$ be a function satisfying Assumptions \ref{assumption_lipschitz} and \ref{assumption_morse}. Further, recall the definition of $\delta$ in \eqref{def_delta} and let $\varepsilon$ be the accuracy of Algorithm \ref{alg_nonconvex_newton}. Then, if \blue{$m<\xi/2$}, the condition $\left\|\nabla f_{-}(\bbx_0) \right\| \geq (5L)/(2m^2)\left\|\nabla f(\bbx_0) \right\|^2$ is violated and $ \left\|\nabla f(\bbx_0) \right\| \leq \zeta \delta$ holds for some $\zeta \in(0,1)$, the iterates either satisfy $\left\|\nabla f(\bbx_{\tilde{k}}) \right\| \leq \varepsilon$ with 
  \begin{equation}
    \tilde{k}=\log_2 \left( \log_2\left(\frac{5L\zeta}{m^2\varepsilon}\right)\right).
  \end{equation}
  or for some $k<\tilde{k}$ we have that $\left\|\nabla f_{-}(\bbx_0) \right\| \leq (5L)/(2m^2)\left\|\nabla f(\bbx_0) \right\|^2$.
\end{proposition}
\begin{proof}
  See Appendix \ref{app_convergence_saddle}.
  \end{proof}
If the iterate $\bbx_0$ is in the $\varepsilon$ neighborhood of a saddle point, then Algorithm 2 {adds uniform noise in each component $X_i\sim \ccalN( 0,2\varepsilon/m)$ to $\bbx_0$. To analyze the perturbed iteration we define the following set. 
\begin{equation}
\ccalG := \left\{\bbx \in \mathbb{R}^n  \big| \left\|\nabla f_{-}(\bbx) \right\| \geq \max\left\{\frac{5L}{2m^2}\left\|\nabla f(\bbx) \right\|^2,\gamma\varepsilon\right\}, \left\|\nabla f(\bbx) \right\| \leq \zeta \delta \right\}.
\end{equation}
We show in Lemma \ref{lemma_probability_good_region} that the probability of $\bbx_0+X\in\ccalG$ is lower bounded by $q$ given in \eqref{eqn_expression_q}. In this case NCN escapes the $\varepsilon$ neighborhood of the saddle at a exponential rate in $1+({\log\left( 1/\gamma\right)})/({\log \left(2-\zeta\right)})$ iterations based on the analysis of Proposition \ref{prop_steps_escape}. 
If the latter is not the case, we show that the number of iterations between re-sample instances is bounded by $O(\log(\log(1/\varepsilon)))$. Because we want to escape each saddle point with probability $(1-p)^{1/S}$ the number of draws needed is of the order of $\log(S/p)$. The following Proposition formalizes the previous discussion.
\begin{proposition}\label{prop_probabilistic}
 Let $f(\bbx)$ be a function satisfying Assumptions \ref{assumption_lipschitz}, \ref{assumption_bounded_opt} and \ref{assumption_morse}. Let \blue{$m<\xi/2$} and $\varepsilon>0$ be the desired accuracy of the output of Algorithm \ref{alg_nonconvex_newton} be such that 
  \begin{equation}
  \varepsilon < \min \left\{\frac{1-\gamma}{4Ln},\frac{\zeta\delta m}{2M\sqrt{n}+m},\frac{2L\gamma}{5\left(2M\sqrt{n}+m\right)^2}\right\}.
  \end{equation}
Further, consider the constants $\gamma, \zeta \in(0,1)$ and let $q$ be a constant such that
\begin{equation}
 q>2\left(1-\Phi(1) \right)\frac{\gamma\left(\frac{n}{2},\frac{n}{2}\right)}{\Gamma\left(\frac{n}{2}\right)},
\end{equation}
  %
  where $\Phi(1)$ is the integral of the Gaussian distribution with integration boundaries $-\infty$ and $1$, $\Gamma(\cdot)$ is the gamma function and $\gamma(\cdot,\cdot)$ is the lower incomplete gamma function. For  any $\bbx_0$ in the neighborhood of a saddle point satisfying $\left\|\nabla f(\bbx_0) \right\| \leq \varepsilon$, with probability $(1-p)^{1/S}$ we have that $\left\|\nabla f_{-}(\bbx_{K_4}) \right\| \geq \frac{5L}{2m^2}\left\|\nabla f(\bbx_{K_4}) \right\|^2$ and $ \varepsilon \leq \left\|\nabla f_{-}(\bbx_{K_4})\right\|$, where $K_4$ is given by 
\begin{equation}
K_4 = \left(1+\frac{\log (S/p)}{ \log \left(1/(1-q)\right)}\right)\left[ \log_2 \left(  \log_2\left(\frac{5L\zeta}{m^2\varepsilon}\right)\right)+\frac{\log \left(1/\gamma\right)}{\log (2-\zeta)}+1\right].
\end{equation} 
\end{proposition}
\begin{proof}
  See Appendix \ref{app_probabilistic}.
\end{proof}
Combining the results from propositions \ref{prop_steps_escape} and \ref{prop_probabilistic} we show that with probability $1-p$ the number of steps required to escape the $\zeta\delta$ neighborhood of a saddle is $O\left(\log(1/p)+\log(1/\varepsilon)\right)$. The previous result completes the analysis of the neighborhoods of the critical points. To complete the convergence analysis we show in the following lemma that the number of iterations required to reach a neighborhood of a critical point satisfying $\left\|\nabla f(\bbx) \right\|<\zeta\delta$ is constant. 
%
\begin{proposition}\label{prop_sublinear}
 Let $f$ be a function satisfying Assumptions \ref{assumption_lipschitz}, \ref{assumption_bounded_opt} and \ref{assumption_morse} and consider $\alpha \in(0,1/2), \beta \in(0,1)$ \blue{and $m<\xi/2$} as the inputs of Algorithm \ref{alg_nonconvex_newton}. Further, recall the definition of $\delta$ in \eqref{def_delta}. Then if $\left\| \nabla f(\bbx_0)\right\|\geq \zeta \delta$ in at most $K_2$ iterations for Algorithm \ref{alg_nonconvex_newton} reaches a neighborhood such that $\left\| \nabla f(\bbx_{K_2})\right\|< \zeta\delta$, with
  \begin{equation}
K_2\leq \frac{M^2\left(f(\bbx_0) -f(\bbx^*)\right)}{\alpha \beta m (\zeta\delta)^2}.
  \end{equation}
\end{proposition}
\begin{proof}
  See Appendix \ref{app_damped_phase}.
  \end{proof}
%
The previous result establishes a bound in the number of iterations that NCN takes to reach a neighborhood of a critical point satisfying $\left\|\nabla f(\bbx) \right\|<\zeta\delta$. However, to complete the proof of the final bound of the complexity of the algorithm, we need to ensure that the algorithm does not visit the neighborhood of the same saddle point indefinitely. In particular the next result estabishes that NCN visits such neighborhoods at most a constant number of times. Moreover, only in one of these visits there is need of adding noise.
%
%
\begin{proposition}\label{prop_revisit_saddle}
  Let $f(\bbx)$ satisfy assumptions \ref{assumption_lipschitz}--\ref{assumption_morse} and consider $\alpha \in(0,1/2), \beta \in(0,1)$ \blue{and $m<\xi/2$} as the inputs of Algorithm \ref{alg_nonconvex_newton}. Let $\bbx_c$ be a critical point of $f(\bbx)$, and let $\zeta \in (0,1)$ and $\delta$ is the constant defined in Theorem \ref{thm_steps_escape} and define 
 \begin{equation}
 \ccalN = \left\{\bbx\in\mathbb{R}^n\Big| \|\bbx-\bbx_c \|\leq m(1-2\alpha)/L, \|\nabla f(\bbx)\| <\zeta \delta\right\}. 
 \end{equation}
 Let $\bbx_0\in\ccalN$ and the desired accuracy $\varepsilon$ satisfy
 \begin{equation}
\varepsilon < \sqrt{\frac{\alpha\beta\frac{m^3}{M^3}(\zeta \delta)^2}{1+\left(2\sqrt{n}\frac{M}{m}+1\right)^2}}.
   \end{equation}
Let $T = \sum_{k=1}^\infty \mathbbm{1}\left(\bbx_k\in\ccalN, \bbx_{k-1}\notin\ccalN\right)$ be the number of times that the sequence generated by NCN visits $\ccalN$. Then, the sequence generated by NCN is such that
 \begin{equation}
 T<\frac{2}{\alpha\beta}\frac{M^3}{m^3}+\alpha\beta,
 \end{equation}
 and the neighborhood $\left\| \nabla f(\bbx)\right\|<\varepsilon$ is visited at most once.\blue{ Likewise, the number of saddle points visited before converging to a local minimum $S$ is upper bounded by
   \begin{equation}
S \leq \frac{m^2}{M\alpha\beta \left(\zeta\delta\right)^2}\left(f(\bbx_0)-f(\bbx^*)\right).
   \end{equation}
   }
\end{proposition}
\begin{proof}
  See Appendix \ref{app_revisit_saddle}.
  \end{proof}
The proof of the final complexity stated in Theorem \ref{main_thm} follows from propositions \ref{prop_steps_escape}, \ref{prop_probabilistic}, \ref{prop_sublinear}, \ref{prop_revisit_saddle}, Corollary \ref{coro_local_minima} and the discussion after the theorem in Section \ref{sec_convergence_preview}. 
%
%
%

{
\section{Numerical Experiments}\label{sec_numerical}
In this section we apply Algorithm \ref{alg_nonconvex_newton} to the matrix factorization problem, where the goal is to find a rank $r$ approximation of a matrix $\bbM \in \ccalM^{l\times n}$, representing ratings from $l=943$ users to $n=1682$ movies \cite{harper2016movielens}. Each user has rated at least $20$ movies for a total of $99900$ known ratings. The problem can be written as
\begin{equation}
(\bbU^*,\bbV^*) := \argmin_{\bbU \in \ccalM^{l\times r}, \bbV \in \ccalM^{n\times r}}  f(\bbU,\bbV)=\argmin_{\bbU \in \ccalM^{l\times r}, \bbV \in \ccalM^{n\times r}}\frac{1}{2} \left\|\bbM-\bbU\bbV^{\top} \right\|_F^2.
\end{equation}
Writing the product $\bbU\bbV^{\top}$ as $\sum_{i=1}^r \bbu_i\bbv_i^{\top}$, where $\bbu_i$ and $\bbv_i$ are the $i$--th columns of the matrices $\bbU$ and $\bbV$ respectively. Let  $\bbx = [\bbu_1^{\top}, \ldots, \bbu_l^{\top}, \bbv_1^{\top}, \ldots \bbv_n^{\top}]^{\top}$, then we solve
\begin{equation}\label{eqn_application}
\bbx^* := \argmin_{\bbx \in \mathbf{R}^{r\times(n+l)}} f(\bbx).
\end{equation}
 We compare the performance of gradient descent and NCN (Algorithm \ref{alg_nonconvex_newton}) on the problem \eqref{eqn_application}. The step size in gradient descent is selected via line search by backtracking with parameters $\alpha$ and $\beta$ being the same as those of NCN. The parameters selected are $r=2$, $\alpha = 0.1$, $\beta=0.9$, accuracy $\varepsilon = 1\times10^{-8}$ and $m=1\times 10^{-12}$. The initial iterate selected for all simulations is the same and it is selected at random from a multivariate normal random variable with mean zero and standard deviation $10$. 

\begin{figure*}
    \centering
    \begin{subfigure}[b]{0.48\textwidth}
        \includegraphics[width=\textwidth]{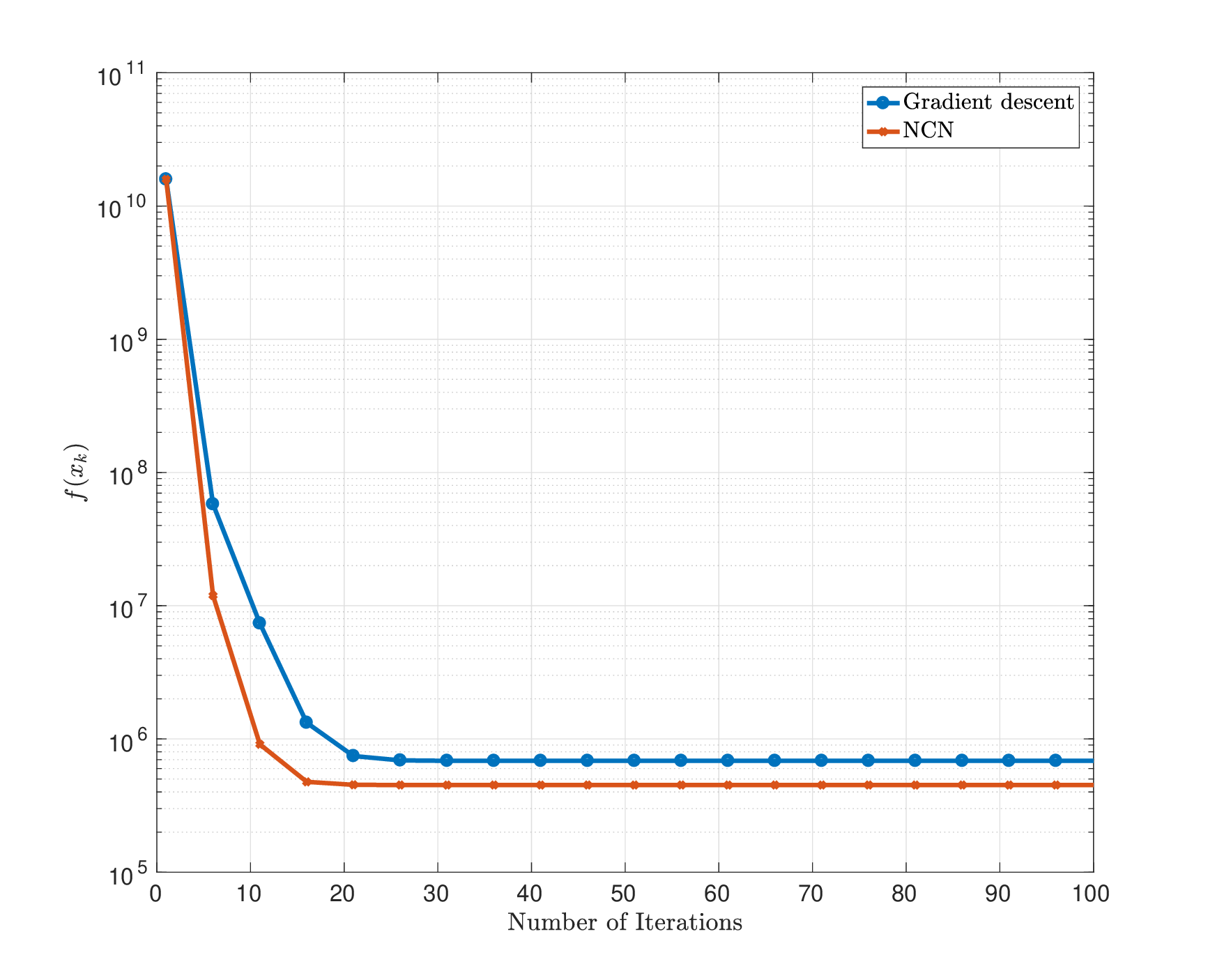}
    \end{subfigure}
    \begin{subfigure}[b]{0.48\textwidth}
        \includegraphics[width=\textwidth]{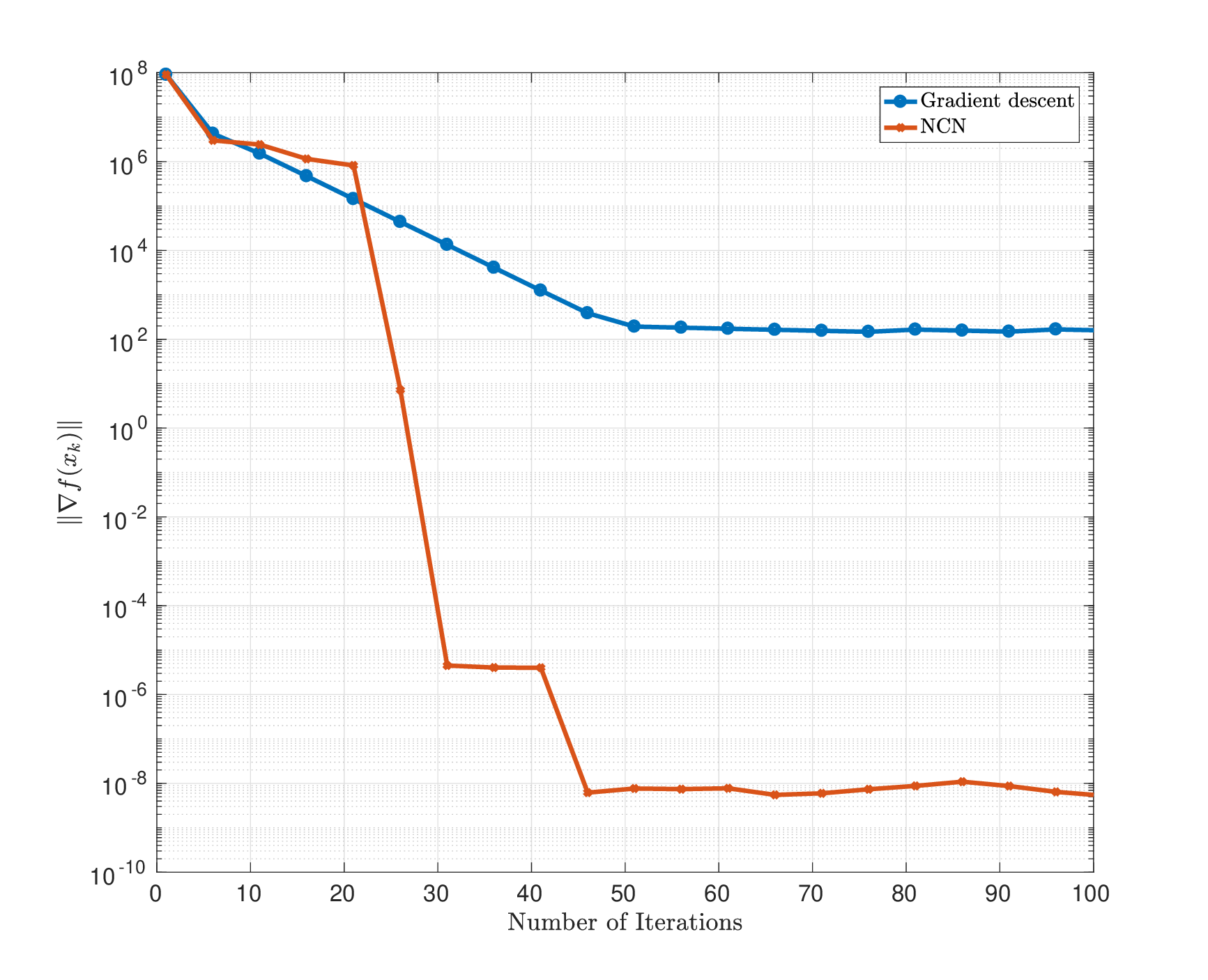}
    \end{subfigure}
    \vspace{-0.4cm}
    \caption{Comparison of gradient descent and NCN with $m = 1\times 10^{-9}$ in terms of value of the objective function (left)  and norm of the gradient (right) for the rank $r= 2$ factorization problem. The dimension of the matrix to factorize is $943$ by $1682$. }\label{fig_movie}
\end{figure*}

  In Figure \ref{fig_movie} we plot the value of the objective function and the norm of the gradient in logarithmic scale for gradient descent and NCN. After $20$ iterations the value of the objective function achieved by NCN is half of that achieved by gradient descent. Notice that the norm of the gradient of the objective function remain large for gradient descent. The latter is an indication of the algorithm being stuck in a neighborhood of a saddle point. Recall from the example in Section \ref{sec_illustrative_example}, that escaping an ill-conditioned saddle point can take an extremely large number of iterations, hence the fact that the norm seems constant is unsurprising. The minimum eigenvalues of the Hessian of the point to which gradient descent converges is $-12.1592$ which indicates that it is converging to a neighborhood of a saddle point. On the other hand, the iterates of Nonconvex Newton achieve a point whose minimum eigenvalue of the Hessian is $-3.0679 \times 10^{-7}$. The latter, even if it is negative and formally the point to which the iterates converge is a saddle point, in practice this point can be considered a local minimum with eigenvalue $0$. Finally, observe that Nonconvex Newton enjoys a fast convergence in the neighborhood of the local minimum, which is to be expected from Corollary \ref{coro_local_minima}.

}

\section{Conclusion}\label{sec_conclusion}
The algorithm presented in this work achieves convergence to a local minimum of a nonconvex function with probability $1-p$ in $O\left(\log(1/p),\log(1/\varepsilon)\right)$ iterations. The main feature of NCN is the curvature correction of the function by pre-multiplying the gradient by the PT-inverse of the Hessian. Since the latter is a positive definite matrix it ensures that the NCN is a descent direction for the objective function and therefore it converges to a neighborhood of the critical points. This feature could be achieved by any other positive definite matrix; however, the structure of PT-inverse enforces a special behavior in the neighborhood of the critical points:  

\begin{itemize}
\item[\bf{(i)}] The projection of the gradient $\nabla f_{+}(\bbx_k)$ exhibits a quadratic convergence behavior. In particular, in the neighborhood of the local minima, this subspace in the whole space, and therefore NCN converges to an $\varepsilon$ neighborhood of the local minimum in $O(\log_2(\log_2(1/\varepsilon)))$ as it is the case for the classic Newton's method.  
\item[\bf{(ii)}] The projection of the gradient over the complement of the previous subspace has an exponential escape rate of $3/2$ independent of the conditioning of the problem as long as a considerable part of the energy of the gradient lies in this subspace. 
\end{itemize}
When the previous condition is not met, it is possible to add noise and with probability $1-p$ after $O\log(1/p)$ iterations the iterate satisfies said condition. The randomization was not needed in the numerical examples, since NCN escapes the saddle before reaching its $\varepsilon$ neighborhood. This suggests that random initialization is enough to ensure an exponential rate of escape in practice, although this is not supported by theoretical guarantees. { The main drawbacks of the proposed algorithm are the need for spectral decomposition -- which is impractical in large scale problems-- and the requirement of saddle points to be non-degenerated -- which might be a strong assumption for many problems of interest. However, the ability to evade saddle points in a number of iterations that is independent of the condition number is a remarkable property of our algorithm, since the condition number of the saddle has an exponential effect in the number of iteration that other algorithms take to escape it. As illusttrated by the example in Section \ref{sec_illustrative_example}. Hence, this is a promising starting point to develop algorithms with better complexity -- computing for instance approximations of the PT-inverse, like quasi-Newton methods do to approximate the Hessian inverse -- and operating under weaker assumptions.

  }


\appendix
\section{Consequences of the Lipschitz continuity of the Hessian}
In this section we state and prove some useful Lemmas needed for the analysis of the convergence of the algorithm.

\begin{lemma}\label{lemma_eigenvalues}
Let $f(\bbx)$ satisfy Assumptions \ref{assumption_lipschitz} and \ref{assumption_morse}. Then the eigenvalues of the Hessian $\lambda_i(\nabla^2 f(\bbx))$ are Lipschitz with constant $L$. 
\end{lemma}
\begin{proof}
We show that for matrices $\bbA, \bbB\in \mathbb{\ccalM}^{n\times n}$, if $\bbA \preceq \bbB$ then $\lambda_i(\bbA)\leq \lambda_i(\bbB)$ for all $i=1\ldots n$. Note that for a given $i$, there exists a subpace $V_1$ of dimension $n-i+1$ where for any $\bbx\in V_1$ we have
  %
$\bbx^\top \bbA\bbx \geq \lambda_i(\bbA)\|\bbx\|^2$.
  %
  Likewise, there exists a subpace $V_2$ of dimension $i$ such that for any $\bbx \in V_2$ we have 
$\bbx^\top \bbB\bbx \leq \lambda_i(\bbB) \|\bbx\|^2$.
  %
  Since the sum of the dimensions of $V_1$ and $V_2$ is larger than $n$. There exists $\bbx \neq \mathbf{0}$, $\bbx\in V_1 \cap V_2$. Since $\bbA \preceq \bbB$ we have that
  %
 $  0\leq \bbx^\top \left(\bbB - \bbA\right) \bbx \leq \left(\lambda_i(\bbB) -\lambda_i(\bbA)\right) \|\bbx\|^2$.
  %
  Therefore $\lambda_i(\bbA)\leq \lambda_i(\bbB)$. Next observe that if $\left\| \bbA-\bbB\right\| \leq \varepsilon $ it holds that $\varepsilon \bbI \pm (\bbA-\bbB)\succeq 0$. Hence, using the previous result we have that
  %
$\varepsilon +\lambda_i(\bbA) \geq \lambda_i(\bbB) \quad \mbox{and} \quad \varepsilon +\lambda_i(\bbB) \geq \lambda_i(\bbA).$
  %
  Which implies that $\left| \lambda_i(\bbB) - \lambda_i(\bbA)\right| \leq \varepsilon$. THe proof is completed by observing that $\nabla^2 f (\bbx)$ is Lipschitz, and therefore for any $\bbx, \bby \in \mathbb{R}^n$ it holds that $\left\|\nabla^2 f (\bbx)-\nabla^2 f (\bby) \right\| \leq L\left\|\bbx-\bby\right\|$. 
  %
  %
  \end{proof}
%
%
\begin{lemma}\label{lemma_bound_eigenvalues}
Let $f(\bbx)$ satisfy Assumptions \ref{assumption_lipschitz} and \ref{assumption_morse} and $\bbx_c$ be a critical point $f(\bbx)$. For any $\bbx \in \mathbb{R}^n$ such that $\left\| \bbx-\bbx_c\right\| <m(1-2\alpha)/L$, with {$m<\zeta/2$} it holds that 
\begin{equation}\label{eqn_eigenvalue_bound}
\min_{i=1\ldots n} \lambda_{i}\left|(\nabla^2 f(\bbx))\right| \geq m \quad \mbox{and} \quad m\left\|\bbx-\bbx_c\right\|\leq \left\| \nabla f(\bbx) \right\|.
\end{equation}
 %
\end{lemma}
\begin{proof}
Let $\bbx_c$ be a critical point with {$|\lambda_{i}(\nabla^2f(\bbx_c))|>\zeta>2m$} (c.f. Assumption \ref{assumption_morse}) and let $\bbx\in \mathbb{R}^n$ be such that $\left\| \bbx - \bbx_c \right\| <m(1-2\alpha)/L$. By the Lipschitz continuity of the eigenvalues (c.f. Lemma \ref{lemma_eigenvalues}), for all $i=1\ldots n$  we have that 
\begin{equation}
\left|\lambda_{i}\left(\nabla^2f(\bbx)\right) - \lambda_{i}\left(\nabla^2 f(\bbx_c)\right) \right| \leq L \left\|\bbx-\bbx_c \right\| \leq (1-2\alpha)m <m.
\end{equation}
Since the difference between the $i$-th eigenvalue at $\bbx$ and $\bbx_c$ is bounded by $m$, the first claim in \eqref{eqn_eigenvalue_bound} holds. 
%
%
To prove the second claim, write the gradient of $f(\bbx)$ as
\begin{equation}
\nabla f(\bbx) = \nabla^2 f(\bbx_c)\left(\bbx-\bbx_c\right)+\frac{1}{2}\left(\nabla^2 f(\tilde{\bbx})-\nabla^2 f(\bbx_c)\right)\left(\bbx-\bbx_c\right),
  \end{equation}
where $\tilde{\bbx} = \theta \bbx_c + (1-\theta)\bbx$ with $\theta \in[0,1]$. Since the minimum absolute value of the eigenvalues of the Hessian at $\bbx_c$ is larger than $2m$ we have that
\begin{equation}\label{eqn_appendixA_need_space}
    \left\|\nabla f(\bbx)\right\| \geq 2m\left\|\bbx-\bbx_c\right\|-\frac{1}{2}\left\|\left(\nabla^2 f(\tilde{\bbx})-\nabla^2 f(\bbx_c)\right)\left(\bbx-\bbx_c\right)\right\|.
  \end{equation}
Using the Lipschitz continuity of the Hessian (c.f. Assumption \ref{assumption_lipschitz}), 
\eqref{eqn_appendixA_need_space} reduces to
  \begin{equation}
  \begin{split}
    \left\|\nabla f(\bbx)\right\| &\geq 2m\left\|\bbx-\bbx_c\right\|-\frac{L}{2}\left\|\bbx-\bbx_c\right\|^2 \\
    &\geq m\left(2-\frac{L\left\|\bbx-\bbx_c\right\|}{2}\right)\left\|\bbx-\bbx_c\right\| \geq m\left\|\bbx-\bbx_c\right\|,
\end{split}
\end{equation}
  where the last inequality follows from the fact that $\left\|\bbx-\bbx_c\right\| <{m}/{L}$.
\end{proof}
\section{Proof of Lemma \ref{lemma_step_size_one}}\label{sec_step_size_one}
Let us define the vector $\Delta \bbx = -\bbH_m(\bbx)^{-1}\nabla f(\bbx)$ and the following function $\tilde{f}:[0,1]\to \mathbb{R}$,
  %
$\tilde{f}(\eta) = f(\bbx+\eta\Delta \bbx)$.
  %
Differentiating with respect to $\eta$ the first two derivatives of $\tilde{f}$ yield
  \begin{equation}
\tilde{f}^\prime(\eta) = \nabla f(\bbx +\eta \Delta \bbx)^T \Delta \bbx \quad
  %
 \mbox{ and }\quad
\tilde{f}^{\prime\prime}(\eta) = \Delta \bbx^T \nabla^2 f(\bbx+ \eta\Delta \bbx)\Delta \bbx.
  \end{equation}
  Then, by virtue of the Lipschitz continuity of $\nabla^2f(\bbx)$ (c.f. Assumption \ref{assumption_lipschitz}) we can upper bound the absolute value of the difference of $\tilde{f}^{\prime \prime}(\eta)$ evaluated at $\eta$ and $0$ by
  \begin{equation}
      |\tilde{f}^{\prime \prime} (\eta) - \tilde{f}^{\prime \prime}(0)| \leq \| \nabla^2f(\bbx+\eta\Delta \bbx)-\nabla^2f(\bbx)\|\|\Delta \bbx\|_2^2 \leq \eta L\|\Delta \bbx\|_2^3.
  \end{equation}
  The latter inequality allows us to upper bound the second derivative of $\tilde{f}(\eta)$ by
  \begin{equation}
\tilde{f}^{\prime \prime}(\eta) \leq \tilde{f}^{\prime \prime}(0) + \eta L\|\Delta \bbx\|_2^3.
  \end{equation}
  Integrating twice with respect to $\eta$ in both sides of the above inequality yields
  \begin{equation}\label{eqn_upper_bound_function_stepsize}
  \tilde{f}(\eta) \leq \tilde{f}(0)+ \tilde{f}^\prime(0)\eta + \tilde{f}^{\prime \prime}(0)\frac{\eta^2}{2} + \frac{\eta^3}{6}L \|\Delta \bbx\|_2^3.
  \end{equation}
Evaluate \eqref{eqn_upper_bound_function_stepsize} at $\eta=1$ to get the following upper bound for $f(\bbx+\Delta \bbx)$
  \begin{equation}\label{eqn_upper_bound_function}
 f(\bbx+\Delta \bbx) \leq f(\bbx) + \nabla f(\bbx)^T \Delta \bbx+\frac{1}{2}\Delta \bbx^T \nabla^2f(\bbx)\Delta \bbx+ \frac{L}{6} \|\Delta \bbx\|_2^3.
  \end{equation}
  We next work towards an upper bound for the term $\Delta \bbx^T \nabla^2f(\bbx)\Delta \bbx$. Using the eigenvalue decomposition it is possible to write this term as
  \begin{equation}
    \Delta \bbx^T \nabla^2f(\bbx)\Delta \bbx =  \left(\bbQ(\bbx)\nabla f(\bbx)\right)^\top \left|\mathbf{\Lambda}(\bbx) \right|_m^{-1} \mathbf{\Lambda}(\bbx) \left|\mathbf{\Lambda}(\bbx) \right|_m^{-1} \bbQ(\bbx)\nabla f (\bbx).
    \end{equation}
  Observe that by definition it holds that $ \mathbf{\Lambda}(\bbx) \preceq \left|\mathbf{\Lambda}(\bbx) \right|_m$, thus we have that
    \begin{equation}\label{eqn_upper_bound_quadratic}
    \Delta \bbx^T \nabla^2f(\bbx)\Delta \bbx \leq \left(\bbQ(\bbx)\nabla f(\bbx)\right)^\top \left|\mathbf{\Lambda}(\bbx) \right|_m^{-1} \bbQ(\bbx)\nabla f (\bbx) = -\nabla f(\bbx)^\top \Delta \bbx.
    \end{equation}
    Next, we derive an upper bound for $\|\Delta \bbx\|_2^3$. By definition of $\Delta \bbx$ we have that
    \begin{equation}
   \|\Delta \bbx\|_2^2 = \nabla f (\bbx)^T \bbH_m(\bbx)^{-1} \bbH_m(\bbx)^{-1} \nabla f(\bbx).
    \end{equation}
    Since the previous expression is maximized when $\nabla f(\bbx) $ is collinear with the eigenvector corresponding to the eigenvalue with minimum absolute value we have that
        \begin{equation}\label{eqn_upper_bound_update}
   \|\Delta \bbx\|_2^2 \leq -\frac{1}{m}\nabla f (\bbx)^T \Delta \bbx.
        \end{equation}
        Combing the upper bounds derived in \eqref{eqn_upper_bound_function}, \eqref{eqn_upper_bound_quadratic} and \eqref{eqn_upper_bound_update} we have that
        \begin{equation}\label{eqn_further_bound_line_search}
f(\bbx+\Delta \bbx) \leq f(\bbx) +\nabla f(\bbx)^T \Delta \bbx\left(\frac{1}{2}-\frac{L}{6m^{3/2}}|-\nabla f(\bbx)^T\Delta \bbx|^{1/2}\right)
        \end{equation}
         We then write the following upper bound for $|-\nabla f(\bbx)^T \Delta \bbx|^{1/2}$ using the fact that the eigenvalue with the minimum absolute value of $\bbH_m(\bbx)$ is $m$
        \begin{equation}
|-\nabla f(\bbx)^T \Delta \bbx|^{1/2}  \leq \frac{1}{m^{1/2}}\|\nabla f(\bbx)\|_2
        \end{equation}
        Which allows us to upper bound \eqref{eqn_further_bound_line_search} by 
\begin{equation}
f(\bbx+\Delta \bbx) \leq f(\bbx)+\frac{1}{2}\nabla f(\bbx)^T \Delta \bbx\left(1-\frac{L}{3m^2}\|\nabla f(\bbx) \|_2\right) 
        \end{equation}
        Finally, for any $\bbx\in\mathbb{R}^n$ such that $\|\nabla f(\bbx)\|_2<({3m^2}/{L})\left(1-2\alpha\right)$ we have that
        \begin{equation}
f(\bbx+\Delta \bbx) \leq f(\bbx) +\alpha \nabla f(\bbx)^T\Delta \bbx,
        \end{equation}
        which shows that line search admits a step size of size $1$.
%
%
%
\section{Proof of Corollary \ref{coro_local_minima}}\label{app_local_minima}
In the neighborhood of a local minimum all the eigenvalues of the Hessian are positive, thus \eqref{eqn_rate_saddles_pos} reduces to 
\begin{equation}
\| \nabla f (\bbx_{k+1})\| =\| \nabla f_{+} (\bbx_{k+1})\| \leq \frac{5L}{2m^2} \left\|\nabla f( \bbx_k) \right\|^2.
\end{equation}
Multiplying both sides of the previous equation by $5L/2m^2$ yields
\begin{equation}
 \frac{5L}{2m^2}\| \nabla f (\bbx_{k+1})\|  \leq \left(\frac{5L}{2m^2} \left\|\nabla f( \bbx_k) \right\|\right)^2.
\end{equation}
The previous inequality can be written recursively as 
\begin{equation}
 \frac{5L}{2m^2}\| \nabla f (\bbx_{k})\|  \leq \left(\frac{5L}{2m^2} \left\|\nabla f( \bbx_0) \right\|\right)^{2^k} \leq 2^{-2^k}, 
\end{equation}
where the right most inequality holds since $\left\|\nabla f(\bbx_0) \right\| <\delta\leq m^2/(5L)$. Thus
\begin{equation}
 \frac{5L}{2m^2}\| \nabla f (\bbx_{K_3})\| \leq 2^{-2^{K_3}} = \frac{5L}{2m^2}\varepsilon. 
\end{equation}
Which completes the proof of the corollary.
%

\section{Proof of Proposition \ref{prop_steps_escape}}\label{app_steps_escape}
Let us start by showing that if $\left\|\nabla f_{-}(\bbx_0) \right\| > (5L)/(2m^2)\left\|\nabla f(\bbx_0) \right\|^2$ then $\left\|\nabla f_{-}(\bbx_1) \right\| > \left\|\nabla f_{+}(\bbx_1) \right\|$. Using the results from Theorem \ref{theo_rate_saddles} we have that
\begin{equation}\label{eqn_bound_x1}
\left\| \nabla f_{-}(\bbx_{1})\right\| \geq 2\left\| \nabla f_{-}(\bbx_{0})\right\|- \frac{5L}{2m^2} \left\| \nabla f(\bbx_{0})\right\|^2 >\left\|\nabla f_{-}(\bbx_{0})\right\|\geq \varepsilon, 
\end{equation}
and
\begin{equation}
\left\| \nabla f_{+}(\bbx_{1})\right\| \leq  \frac{5L}{2m^2} \left\| \nabla f(\bbx_{0})\right\|^2 <\left\|\nabla f_{-}(\bbx_{0})\right\|.  
\end{equation}
Therefore $\left\|\nabla f_{-}(\bbx_1) \right\| > \left\|\nabla f_{+}(\bbx_1) \right\|$. We next next show that the same holds for all $k\geq 1$. Let us prove it by induction. Assume $\left\| \nabla f_{-}(\bbx_{k})\right\| \geq \left\| \nabla f_{+}(\bbx_{k})\right\|$ holds for some $k\geq 1$. In this case we have that  $\left\|\nabla f(\bbx_k) \right\| \leq 2\left\|\nabla f_{-}(\bbx_k) \right\|$, thus
\begin{equation}\label{eqn_bound_for_constant_rate_new}
\frac{5L}{2m^2} \left\| \nabla f(\bbx_{k})\right\|^2 \leq \frac{5L}{m^2} \left\| \nabla f(\bbx_{k})\right\|\left\| \nabla f_{-}(\bbx_{k})\right\| \leq \zeta \left\| \nabla f_{-}(\bbx_{k})\right\|,
\end{equation}
where, the last inequality follows from the fact that $\left\|\nabla f(\bbx_k) \right\| \leq \zeta \delta \leq \zeta (m^2)/(5L)$. Since $\left\|\nabla f(\bbx_k) \right\| \leq \zeta \delta$ the hypothesis of Theorem \ref{theo_rate_saddles} are satisfied and we have that 
\begin{equation}\label{eqn_negative_constant_rate}
\left\| \nabla f_{-}(\bbx_{k+1})\right\| \geq 2\left\| \nabla f_{-}(\bbx_{k})\right\|- \frac{5L}{2m^2} \left\| \nabla f(\bbx_{k})\right\|^2 \geq \left(2-\zeta\right) \left\| \nabla f_{-}(\bbx_{k})\right\|
\end{equation}
and 
\begin{equation}\label{eqn_positive_constant_rate}
\left\| \nabla f_{+}(\bbx_{k+1})\right\| \leq \frac{5L}{2m^2} \left\|\nabla f(\bbx_k) \right\|^2\leq \zeta \left\| \nabla f_{-}(\bbx_{k})\right\|,
\end{equation}
where the rightmost inequality in the two previous expressions follows from the bound \eqref{eqn_bound_for_constant_rate_new}. The latter implies that $\left\| \nabla f_{+}(\bbx_{k+1})\right\| \leq \left\| \nabla f_{-}(\bbx_{k+1})\right\|$ since $\zeta\in(0,1)$. Which proves that for all $k\geq 1$ we have that $\left\| \nabla f_{+}(\bbx_{k})\right\| \leq \left\| \nabla f_{-}(\bbx_{k})\right\|$. The latter allows us to write \eqref{eqn_negative_constant_rate} for every $k\geq 1$ as long as $\left\| \nabla f(\bbx_k) \right\| \leq \zeta \delta$. 
Therefore, writing \eqref{eqn_negative_constant_rate} recursively, for $K_2>0$ it holds that
\begin{equation}
\left\| \nabla f_{-}(\bbx_{K_2})\right\| \geq  \left(2-\zeta\right)^{K_2-1} \left\| \nabla f_{-}(\bbx_{1})\right\| \geq \left(2-\zeta\right)^{K_2-1} \varepsilon,
\end{equation}
where we used \eqref{eqn_bound_x1} to lower bound $\left\|\nabla f_{-}(\bbx_1)\right\|>\varepsilon$ in the rightmost inequality in the previous expression. Hence, for $\bbx_{K_1}$, with $K_2 =1+ \log\left( \zeta \delta/\varepsilon\right) /\log \left(2-\zeta\right)$ we have 
\begin{equation}
\left\| \nabla f_{-}(\bbx_{K_2})\right\| \geq \left(2-\zeta\right)^{\frac{\log\left( \zeta \delta/\varepsilon\right)}{\log \left(2-\zeta\right)}} \varepsilon = e^{\log \frac{\zeta \delta}{\varepsilon}}\varepsilon = \zeta \delta.
\end{equation}
Which completes the proof of the Proposition.

 %
%
\section{Proof of Proposition \ref{prop_quadratic_convergence}}\label{app_convergence_saddle}
Consider the case in which $\left\|\nabla f_{-}(\bbx_k) \right\|< \frac{5L}{2m^2}\left\| \nabla f(\bbx_k)\right\|^2$ for all $k$. Because $\left\|\nabla f(\bbx_k)\right\| \leq \delta\leq m^2/(5L)$ we have that
\begin{equation}
\left\|\nabla f_{-}(\bbx_k) \right\|< \frac{1}{2}\left\| \nabla f(\bbx_k)\right\|.
\end{equation}
Hence, $\left\|\nabla f_{-}(\bbx_k) \right\|< \left\|\nabla f_{+}(\bbx_k) \right\|$.  This, being the case for all $k$ implies that
\begin{equation}
\left\| \nabla f(\bbx_{k+1})\right\|\leq 2\left\| \nabla f_{+}(\bbx_{k+1})\right\| \leq \frac{5L}{m^2}\left\| \nabla f(\bbx_{k})\right\|^2,
\end{equation}
where the last inequality follows from the result of Theorem \ref{theo_rate_saddles}. Then write recursively the previous expression as
\begin{equation}
\frac{5L}{m^2}\left\| \nabla f(\bbx_{k})\right\|\leq \left(\frac{5L}{m^2}\left\| \nabla f(\bbx_{k})\right\|\right)^{2^k} \leq \zeta^{2^k}.
\end{equation}
Thus, in at most $\tilde{k} = \log_2 \left( \log_2\left(\frac{5L\zeta}{m^2\varepsilon}\right)\right)$ we have that $\left\| \nabla f(\bbx_{k+\tilde{k}}) \right\| \leq \varepsilon$.
%
%
\section{Proof of Proposition \ref{prop_probabilistic}}\label{app_probabilistic}
Let $\bbx_k$ be such that $\left\| \nabla f(\bbx_k)\right\| < \varepsilon$ and let $q$ be the probability of having $\bbx_k +X \in \ccalG$, with $X$ a Gaussian random vector where each component is given by $X_i \sim \ccalN(0, 2\varepsilon/m)$. If this is the case, the rate of escape of the $\varepsilon$ neighborhood is $2-\zeta$ and therefore we escape this region in $1+\log(1/\gamma)/\log(2-\zeta)$ iterations. The second case is when $\bbx_k +X \notin \ccalG$, which happens with probability $1-q$. Two possibilities arise in this scenario, either $\left\| \nabla f(\bbx_k+X)\right\|<\varepsilon$ or not. If the former happens then Algorithm \ref{alg_nonconvex_newton} will perform $1+\log(1/\gamma)/\log(2-\zeta)$ iterations and it might be case that the argument is still in the $\varepsilon$ neighborhood. Thus requiring to perturb the iterate again as stated by Algorithm \ref{alg_nonconvex_newton}. The other possibility is that either immediately after the perturbation or after the $1+\log(1/\gamma)/\log(2-\zeta)$ iterations we have $\left\| \nabla f(\bbx_k)\right\| >\varepsilon$, but $\bbx_k$ is still not in $\ccalG$. In this situation from Proposition \ref{prop_quadratic_convergence} it yields that in at most $\tilde{k}$ iterations we have that $\left\| \nabla f(\bbx_{\tilde{k}})\right\| <\varepsilon$, where $\tilde{k}$ is given by 
\begin{equation}
\tilde{k} =\log_2 \left(  \log_2\left(\frac{5L\zeta}{m^2\varepsilon}\right)\right).
\end{equation}
 To summarize, with probability $1-q$ the perturbation $X$ yields an iterate such that $\bbx_k+X \notin\ccalG$. This being the case, after at most $k^*$ iterations the argument is in the $\varepsilon$ neighborhood and thus we draw a new random variable, with 
\begin{equation}
k^*= \log_2 \left(  \log_2\left(\frac{5L\zeta}{m^2\varepsilon}\right)\right)+\frac{\log \left(1/\gamma\right)}{\log (2-\zeta)}+1.
\end{equation}
We are left to show that the number of draws needed to escape each saddle with probability $(1-p)^{1/S}$ is given by $\log(S/p)/\log(1/(1-q))$.  Since the draws are independent, the probability of having $x_k \notin\ccalG$ after $N$ perturbations is given by $(1-q)^N$. Hence, we require
%
$(1-q)^N < 1-(1-p)^{\frac{1}{S}} \approx {p}/{S}$ draws,
%
where we have used the Taylor expansion of $(1-p)^{1/S}$ around the point $p=0$. Taking logarithm on both sides of the previous inequality yields
\begin{equation}
N> \frac{\log(p/S)}{\log(1-q)} = \frac{\log(S/p)}{\log(1/(1-q))}.
\end{equation}
Thus completing the proof of the proposition. 
\begin{lemma}\label{lemma_probability_good_region}
    Let $f(\bbx)$ be a function satisfying Assumptions \ref{assumption_lipschitz} and \ref{assumption_morse}. Let $X$ be a random vector, such that $X_i \sim \ccalN(0,2\varepsilon/m)$ and $\bbx_0 \in \mathbb{R}^n$ in the neighborhood of a saddle point such that $\left \|\nabla f(\bbx_0)\right\| < \varepsilon$ with 
\begin{equation}
\varepsilon < \min \left\{\frac{1-\gamma}{4Ln},\frac{\zeta\delta m}{2M\sqrt{n}+m},\frac{2L\gamma}{5\left(2M\sqrt{n}+m\right)^2}\right\}.
\end{equation}
    Then $q=P\left(\bbx_0+X\in \ccalG\right)$ satisfies \eqref{eqn_expression_q}.
    %
%
\end{lemma}
\begin{proof}
Let us define the following set for $\alpha>1$
\begin{equation}
\ccalG_{\alpha} := \left\{\bbx \in\mathbb{R}^n \Big| \left\|\nabla f_{-}(\bbx)\right\|>\gamma \varepsilon, \left\|\nabla f(\bbx)\right\|<\alpha\gamma \varepsilon,\left\|\nabla f_{-}(\bbx)\right\|\geq \frac{5m^2}{2L}\left\|\nabla f(\bbx)\right\|^2 \right\},
\end{equation}
and the following value for $\alpha_1$
\begin{equation}
\alpha_1 = \sqrt{\frac{2L}{5m^2\gamma\varepsilon}} > \frac{2M\sqrt{n}+m}{m\gamma}>\frac{1}{\gamma}.
\end{equation}
Note as well that for every $\alpha <\alpha_1$ if $\left\|\nabla f_{-}(\bbx)\right\|>\gamma\varepsilon$ and $\left\|\nabla f(\bbx)\right\|<\alpha\gamma \varepsilon$, then $\left\|\nabla f_{-}(\bbx)\right\|\geq \frac{5m^2}{2L}\left\|\nabla f(\bbx)\right\|^2$ is satisfied. Hence, for any $\alpha <\alpha_1$ we have that
\begin{equation}
\ccalG_{\alpha} := \left\{\bbx \in\mathbb{R}^n \Big| \left\|\nabla f_{-}(\bbx)\right\|>\gamma \varepsilon, \left\|\nabla f(\bbx)\right\|<\alpha\gamma \varepsilon \right\}.
\end{equation}
In addition if $\alpha<\alpha_2= (\zeta\delta)/(\gamma\varepsilon) $, $\ccalG_\alpha\subset \ccalG$. Hence for any $\alpha \leq \alpha_2$ it holds that
%
$P\left(\bbx_0+X\in\ccalG\right)>P\left(\bbx_0+X\in\ccalG_\alpha\right).$
%
and, because $\varepsilon < (\zeta\delta m)/(5M\sqrt{n}+1)$,
\begin{equation}
\alpha_2>\frac{2M\sqrt{n}+m}{m\gamma}  >\frac{1}{\gamma},
\end{equation}
 For what remains we chose $\alpha^* = (2M\sqrt{n}+m)/(m\gamma).$ 
We proceed to lower bound $\|\nabla f_{-}(\bbx_0+X)\|$. From its definition \eqref{def_grad_neg_and_pos} and the Taylor expansion around $\bbx_c$ we have
\begin{equation}
\nabla f_{-}(\bbx_0+X) = \bbQ_{-}^\top \nabla f(\bbx_0+X) = \bbQ_{-}^\top \nabla^2f(\tilde{\bbx})(\bbx_0+X-\bbx_c),
\end{equation}
where $\tilde{\bbx} = \theta\bbx+(1-\theta)\bbx_c$ for some $\theta \in [0,1]$. Using the Lipschitz continuity of the Hessian and we can upper bound the norm of the previous expression by
\begin{equation}
\left\|\nabla f_{-} (\bbx_0+X) \right\| \geq \left\| \bbQ_{-}^\top \nabla^2 f(\bbx_c)\left(\bbx_0+X-\bbx_c\right)\right\| - L \left\|\bbx_0+X-\bbx_c\right\|^2.
\end{equation}
Since the absolute value of the minimum eigenvalue of the Hessian at $\bbx_c$ is at least $2m$ (c.f. Assumption \ref{assumption_morse}) it holds that
\begin{equation}\label{eqn_bound_nabla_neg}
\begin{split}
\left\|\nabla f_{-} (\bbx_0+X) \right\| &\geq 2m\left\| \bbQ_{-}^\top \left(\bbx_0+X-\bbx_c\right)\right\| - L \left\|\bbx_0+X-\bbx_c\right\|^2\\
& \geq 2m\left\| \bbQ_{-}^\top X\right\| -2m\left\| \bbx_0-\bbx_c\right\| - L \left\|\bbx_0-\bbx_c\right\|^2-L\left\|X\right\|^2.
\end{split}
\end{equation}
Recall that (c.f. Section \ref{sec_step_size_one}) if $\left\| \nabla f(\bbx_0)\right\| < \delta $ we have that 
%
$m \left\| \bbx_0-\bbx_c\right\| \leq \left\| \nabla f(\bbx_0) \right\| \leq {m^2}/{L}$,
%
where the last inequality follows from the definition of $\delta$ in \eqref{def_delta}. Hence
\begin{equation}
2m\left\|\bbx_0-\bbx_c\right\|+L\left\|\bbx_0-\bbx_c \right\|^2 \leq 3m \left\|\bbx_0-\bbx_c\right\|\leq 3 \left\| \nabla f(\bbx_0)\right\| < 3\varepsilon.
\end{equation}
Combining the previous bound with \eqref{eqn_bound_nabla_neg} yields
\begin{equation}
\left\|\nabla f_{-} (\bbx_0+X) \right\| \geq 2m\left\| \bbQ_{-}^\top X\right\|-L\left\|X\right\|^2-3\varepsilon.
\end{equation}
On the other hand, the Lipschitz continuity of the gradient (c.f. Assumption \ref{assumption_lipschitz}) implies
\begin{equation}
\left\| \nabla f(\bbx_0+X)\right\| \leq \left\| \nabla f(\bbx_0)\right\|+ M\left\|X\right\| \leq \varepsilon + M\left\|X\right\|.
\end{equation}
Then, define $\beta =(\alpha^*\gamma-1)\varepsilon/M$.  Note that $\beta>0$ since $\alpha>1/\gamma$ and write  
\begin{equation}
\begin{split}
P\left(\bbx_0+X\in\ccalG_{\alpha^*}\right) >P\left(\left\|\bbQ_{-}^\top X\right\| -\frac{L}{2m}\|X\|^2> (3+\gamma)\frac{\varepsilon}{2m}, \left\|X\right\|<\beta\right),
\end{split}
\end{equation}
Conditioning on $\left\|X\right\|<{\beta}$ we can further lower bound the previous probability by
\begin{equation}\label{eqn_probability_first_bound}
P\left(\bbx_0+X\in\ccalG_{\alpha^*}\right) > P\left( \left\|X\right\|<{\beta}\right)\ P\left(\left\| \bbQ_{-}^\top X\right\| > (3+\gamma)\frac{\varepsilon}{2m} +\frac{L}{2m}{\beta}^2\right).
\end{equation}
Let us next derive a lower bound for $P(\left\|X\right\|<{\beta})$. Since $X_i$ are normal with mean zero and standard deviation $\sigma$, $\|X\|^2/\sigma^2$ is Chi-squared with parameter $n$. Thus, 
\begin{align}
P\left(\left\|X\right\|<{\beta}\right)  = P\left(\frac{\|X\|^2}{\sigma^2} < \frac{{\beta}^2}{\sigma^2}\right) 
&= P\left(\frac{\|X\|^2}{\sigma^2} < \frac{(\alpha^*\gamma-1)^2\varepsilon^2}{M^2\sigma^2}\right) 
\nonumber\\
&=P\left(\frac{\|X\|^2}{\sigma^2} < \frac{(\alpha^*\gamma-1)^2m^2}{4M^2}\right),
\end{align}
where in the last equality we have used the fact that $\sigma:=2\varepsilon/m$. Since $\alpha^*=(2M\sqrt{n}+m)/(\gamma m)$ we have that $(\alpha^*\gamma-1)m = 2M\sqrt{n}$. Hence
\begin{equation}\label{eqn_bound_norm_proba}
P\left(\left\|X\right\|<{\beta}\right)  > P\left(\frac{\|X\|^2}{\sigma^2}<n\right)=\frac{\gamma\left(\frac{n}{2},\frac{n}{2}\right)}{\Gamma\left(\frac{n}{2}\right)}.
\end{equation}
Note that the ratio ${\gamma\left(\frac{n}{2},\frac{n}{2}\right)}/{\Gamma\left(\frac{n}{2}\right)}$ for any $n\geq1$ is lower bounded by $1/2$. Therefore, we obtain that $P\left(\left\|X\right\|<{\beta}\right)>1/2$. Now we proceed to derive a lower bound for the second probability in the right hand side of \eqref{eqn_probability_nabla_neg}. Let $l$ be the number of negative eigenvalues of the Hessian and $\bbq_i$ be the eigenvectors corresponding to those eigenvalues and write the norm $\left\|\bbQ_{-}^\top X\right\|$ can be written as
\begin{equation}
\left\|\bbQ_{-}^\top X\right\| = \sqrt{\sum_{i=1}^l |\bbq_i^\top X|^2}.
\end{equation}
For each $i$ the product $\bbq_i^\top X$ is a linear combination of Gaussian random variables. Thus, $Y_i = \bbq_i^\top X \sim \ccalN\left(\sum_{j=1}^n \bbq_{ij}\mu_j, \sum_{j=1}^n \bbq_{ij}^2\sigma_j^2\right)$. Since each $X_i$ is normal with mean zero and variance $\sigma^2$ and the norm of $\bbq_i$ is one we have that $Y_i \sim \ccalN(0,\sigma^2)$. Because the minimum norm of $\bbQ_{-}^\top X$ is attained when there is only one negative eigenvalue it is possible to lower bound its norm by$\left\|\bbQ_{-}^\top X \right\| \geq |Y_i|$. Thus we have that
\begin{equation}
P\left(\left\| \bbQ_{-}^\top X\right\| > (3+\gamma)\frac{\varepsilon}{2m} +\frac{L}{2m}{\beta}^2\right)  = P\left(\frac{|Y_i|}{\sigma}>(3+\gamma)\frac{\varepsilon}{2m\sigma}+\frac{L}{2m\sigma}\beta^2\right)
\end{equation}
Recall that $\sigma=2\varepsilon/m$, thus
\begin{equation}
(3+\gamma)\frac{\varepsilon}{2m\sigma } +\frac{L\beta^2}{2m\sigma} = \frac{3+\gamma}{4} +\frac{L(\alpha^*\gamma-1)^2\varepsilon}{4M^2} = \frac{3+\gamma}{4} + Ln\varepsilon<1.
\end{equation}
Where the latter inequality comes from the fact that $\varepsilon<(1-\gamma)/(4Ln)$. Therefore,
\begin{equation}\label{eqn_probability_nabla_neg}
P\left(\left\| \bbQ_{-}^\top X\right\| > (3+\gamma)\frac{\varepsilon}{2m} +\frac{L}{2m}{\beta}^2\right)  \geq 2\left(1-\Phi(1)\right).
\end{equation}
Hence combining the bounds in \eqref{eqn_probability_first_bound}, \eqref{eqn_bound_norm_proba} and  \eqref{eqn_probability_nabla_neg} we have that
\begin{equation}
P(\bbx_0+X\in\ccalG)> P(\bbx_0+X\in\ccalG_{\alpha^*})>2\left(1-\Phi(1)\right)\frac{\gamma\left(\frac{n}{2},\frac{n}{2}\right)}{\Gamma\left(\frac{n}{2}\right)}.
\end{equation}
Thus, completing the proof of the lemma.
\end{proof}

%
\section{Proof of Proposition \ref{prop_sublinear}}\label{app_damped_phase}
Use the fact that the maximum eigenvalue of the Hessian is bounded by $M$ (c.f. Assumption \ref{assumption_lipschitz}) to write
  \begin{equation}
f(\bbx_{k+1}) \leq f(\bbx_k)+\eta_k\nabla f(\bbx_k)^\top\Delta \bbx +\frac{\eta_k^2}{2}M\left\| \Delta \bbx\right\|^2_2.
    \end{equation}
  The previous expression can be further upper bounded using \eqref{eqn_upper_bound_update} as
  \begin{equation}
    f(\bbx_{k+1}) \leq f(\bbx_k)+\eta_k\nabla f(\bbx_k)^\top\Delta \bbx\left(1-\frac{\eta_k M}{2m}\right) 
  \end{equation}
  Observe that $\hat{\eta}=m/M$ satisfies the exit condition of the line search. Hence, backtracking outputs a step size satisfying
  %
$\eta_k>\beta{m}/{M}.$
  %
Thus, we have that
  \begin{equation}\label{eqn_gradient_convergence}
    f(\bbx_{k+1}) \leq  f(\bbx_k)+\alpha\eta_k\nabla f(\bbx_k)^\top\Delta \bbx \leq f(\bbx_k)+\alpha\beta\frac{m}{M}\nabla f(\bbx_k)^\top\Delta \bbx.
  \end{equation}
Observe that maximum eigenvalue of $\nabla^2 f(\bbx)$ is bounded by $M$ hence we have that $-\nabla f(\bbx_k)^\top\Delta \bbx\geq \left\|\nabla f(\bbx_k)\right\|^2/M$. Combine this fact with the characterization of the region away of the saddle points, i.e., $\left\| \nabla f(\bbx_k)\right\|> \delta \zeta$, to bound \eqref{eqn_gradient_convergence} as
  \begin{equation}\label{eqn_damped_decrement}
    f(\bbx_{k+1}) \leq  f(\bbx_k)-\alpha\beta \frac{m}{M^2}(\zeta\delta)^2 .
  \end{equation}
  The previous inequality shows that $f(\bbx)$ is decreasing along the sequence $\left\{\bbx_k\right\}$. Since $f(\bbx)$ is lower bounded the sequence $f(\bbx_k)$ converges, which by virtue of \eqref{eqn_gradient_convergence} implies that $\lim_{k\to\infty}\left\|\nabla f(\bbx_k)\right\| =0$. Let $\bbx^* := \argmin_{\bbx\in\mathbb{R}^n}  f(\bbx)$ and $\bbx_{k}$ be the first iterate satisfying $\left\|\nabla f(\bbx_{k})\right\| \leq \zeta\delta$. Then we can write 
  \begin{equation}
  \begin{split}
f(\bbx_0) -f(\bbx^*)&\geq f(\bbx_0)-f(\bbx_k)= \sum_{i=0}^{k-1} f(\bbx_i)-f(\bbx_{i+1}) \\
& \geq \sum_{i=0}^{k-1} \alpha\beta \frac{m}{M^2}(\zeta\delta)^2  = k \alpha\beta \frac{m}{M^2}(\zeta\delta)^2.
\end{split}
  \end{equation}
Where the last inequality follows from \eqref{eqn_damped_decrement}.  Hence the number of iterations needed for Algorithm \ref{alg_nonconvex_newton} to reach the neighborhood defined by $\left\| \nabla f(\bbx_k)\right\|<\zeta\delta$ is bounded by
  \begin{equation}
k\leq \frac{M^2\left(f(\bbx_0) -f(\bbx^*)\right)}{\alpha \beta m (\zeta\delta)^2}.
  \end{equation}
  \section{Proof of Proposition \ref{prop_revisit_saddle}}\label{app_revisit_saddle}
  Let us start by showing that if we visit a saddle point multiple times, only in one of those visits NCN adds noise to the iterates. Observe that in order to perform a random draw it must hold that $\left\| \nabla f(\bbx)\right\| <\varepsilon$ (c.f. Step 4 of Algorithm \ref{alg_nonconvex_newton}). Using the fact that the maximum eigenvalue of $\nabla^2 f(\bbx)$ is bounded by $M$ we can lower bound $f(\bbx)$ in that neighborhood by
  \begin{equation}
f(\bbx)\geq f(\bbx_c)-\frac{M}{m^2}\varepsilon^2.
  \end{equation}
  Hence, if we show that after adding noise and escaping the saddle, the value of the function is smaller than the right hand side of the above equation, since the update ensures decrement,  the iterates can never reach the neighborhood $\left\| \nabla f(\bbx)\right\| <\varepsilon$ again. Therefore noise is added only once. Observe that, since the random noise is added to ensure that  $\left\| \nabla f(\bbx)\right\|$ is bounded by $(2\sqrt{n}\frac{M}{m}+1)\varepsilon$ (c.f. Step 6 of Algorithm \ref{alg_nonconvex_newton}), it is possible to upper bound $f(\bbx_k)$, with $\bbx_k$ being the iterate after a random draw, by
  \begin{equation}
    f(\bbx_k) \leq f(\bbx_c)+\frac{M}{m^2}\left(2\sqrt{n}\frac{M}{m}+1\right)^2\varepsilon^2. 
    \end{equation}
   Since the line search routine ensures that the function is decreasing along the sequence of iterates generated by NCN (c.f. Step 4 of Algorithm \ref{alg_ncn}) the value of the function of the first iterate outside the neighborhood $\ccalN$ is at most of $f(\bbx_c)+M/m^2\left(2\sqrt{n}\frac{M}{m}+1\right)^2\varepsilon^2$. As in the proof of Proposition \ref{prop_sublinear} we have that the decrement in the value of $f(\bbx_k)$ in the region satisfying $\left\|\nabla f(\bbx)\right\| >\zeta \delta$ is lower bounded by $\alpha \beta (m/M^2)(\zeta\delta)^2$. Hence, if we return to $\ccalN$ in one iteration, it must be the case that
  \begin{equation}
    f(\bbx_k) \leq f(\bbx_c)+\frac{M}{m^2}\left(2\sqrt{n}\frac{M}{m}+1\right)^2\varepsilon^2-\alpha \beta \frac{m}{M^2}(\zeta\delta)^2 < f(\bbx_c)-\frac{M}{m^2}\varepsilon^2.
    \end{equation}
  Where the last inequality follows from the assumption on $\varepsilon$. Therefore the iterates of NCN never reach the neighborhood $\left\| \nabla f(\bbx)\right\| <\varepsilon$ again hence ensuring that only one draw is performed. Using the fact that the function $f$ is decreasing along the sequence of iterates generated by NCN it is possible to show that one can visit each neighborhood of a saddle point only a finite number of times. Let us upper and lower bound the value of $f(\bbx)$ with $\bbx \in \ccalN$ by
  \begin{equation}
f(\bbx_c)-\frac{M}{m^2}(\zeta\delta)^2 \leq f(\bbx) \leq f(\bbx_c)+\frac{M}{m^2}(\zeta\delta)^2.
    \end{equation}
  To make sure that the iterates escape $\ccalN$ we need to ensure that the value of the function along the sequence of iterates generated by NCN decreases its value by $2M/m^2(\zeta\delta)^2 + M/m^2\left(2\sqrt{n}\frac{M}{m}+1\right)^2\varepsilon^2$. The worst case scenario is when after one iteration outside the region $\ccalN$ we are back in it. In each one of these iterations the function value is reduced by $\alpha\beta m/M (\zeta \delta)^2$. Hence the maximum number of times that one visits the same neighborhood of a critical point is upper bounded by
  \begin{equation}
T<\frac{M^3}{m^3}\left(\frac{2}{\alpha\beta} + \left(2\sqrt{n}\frac{M}{m}+1\right)^2\frac{\varepsilon^2}{\left(\zeta \delta\right)^2}\right) < \frac{2}{\alpha\beta}\frac{M^3}{m^3}+\alpha\beta.
    \end{equation}
  {By virtue of the same argument, two consecutive visited critical points satisfy
\begin{equation}
f(\bbx_{c_i}) - f(\bbx_{c_{i+1}}) > \alpha\beta\frac{M}{m^2} \left(\zeta\delta\right)^2.
\end{equation}
Thus, the maximum number of saddle points that can be visisted $S$ satisfies
\begin{equation}
f(\bbx_0)-f(\bbx^*) > S\min_{i=1\ldots S}\left(f(\bbx_{c_i})-f(\bbx_{c_{i+1}})\right) > S\alpha\beta\frac{M}{m^2}\left(\zeta\delta\right)^2, 
\end{equation}
where $\bbx_{c_{S+1}}$ is a local minima. 
}
Thus completing the proof of the proposition. 
\bibliographystyle{siamplain}
\bibliography{bib}



\end{document}